\documentclass{conm-p-l}
\usepackage{amsmath, amssymb, amsthm, graphics}

\theoremstyle{plain}
\newtheorem{prop}{Proposition}[section]
\newtheorem{coro}[prop]{Corollary}
\newtheorem{theo}[prop]{Theorem}
\newtheorem{lemm}[prop]{Lemma}

\theoremstyle{definition}
\newtheorem{defi}[prop]{Definition}

\newtheorem{exam}[prop]{Example}
\newtheorem{ques}[prop]{Question}
\newtheorem{rema}[prop]{Remark}

\numberwithin{equation}{section}
\numberwithin{figure}{section}

\def\Reff#1; #2; #3; #4; #5; #6; #7\par{%
\bibitem{#1} #2, {\it #3}, #4 {\bf #5} (#6) #7}

\def\Ref#1; #2; #3; #4\par{%
\bibitem{#1} #2, {\it #3}, #4}

\def\a{a}
\def\AM{A_M}

\def\b{b}

\def\c{c}
\def\card{{\rm card}}
\def\CC{C}
\def\cov(#1){\mathrel{\triangleright\!
\raise-3pt\hbox{$\scriptstyle #1$}}}

\let\D=\Delta
\def\dist{{\rm dist}}
\def\dive{\preceq}
\def\div{\prec}
\def\Div{{\rm Div}}
\def\Dl{{\rm Div}}
\def\Dr{{\rm Div}_r}
\def\dr{\mathord{\backslash}}

\def\G{\Gamma}
\def\Gr(#1){\langle #1\rangle}

\def\ie{{\it i.e.}}
\def\ii{^{-1}}
\def\ince{\subseteq}

\def\Mon(#1){\langle #1\rangle^{\scriptscriptstyle\!+}}

\def\n(#1){\Vert #1 \Vert}
\def\NN{{\mathbb N}}
\def\notcov(#1){\mathrel{\not\triangleright\!
\raise-3pt\hbox{$\scriptstyle #1$}}}

\def\op{}

\def\PM{P_M}
\let\pp=\dots

\def\PresI{\a, \b \, ; \,
   \a^2 = \b^2, \a\b = \b\a}
\def\PresII{\a, \b, \c \, ; \, 
   \a^2 = \b^2 = \c^2,\a\b = \b\c = \c\a, \a\c = \b\a =
    \c\b}
\def\PresIII{\a, \b, \c \, ; \,
    \a\c = \c\a = \b^2, \a\b = \b\c, \c\b =  \b\a}
    
\def\q{quasi}

\def\resp{\hbox{resp.} }
\def\RRX{R_\SS^*}
\def\RX{R_\SS}

\let\s=\sigma
\let\S=\Sigma  
\def\SD{\S_\D}
\def\ss{s}
\def\SS{S}  
\def\SSU{\SS\op\UM}

\def\TT{T}  

\def\uu{u}

\def\vv{v}

\def\ww{w}

\def\xx{x}

\def\yy{y}
\def\YY{Y}

\def\zz{z}
\def\ZZ{{\mathbb Z}}

\def\dd{^*}
\def\ddd{^{**}}
\def\UM{M^*}
\def\q{quasi}
\def\Q{Quasi}

\begin{document}

\author{Patrick DEHORNOY}
\address{Laboratoire SDAD, Math\'ematiques\\
Universit\'e de Caen BP 5186, 14032 Caen, France}
\email{dehornoy@math.unicaen.fr}
\urladdr{//www.math.unicaen.fr/\!\hbox{$\sim$}dehornoy}

\title{THIN GROUPS OF FRACTIONS}



\begin{abstract}
  A number of properties of spherical Artin groups
  extend to Garside groups, defined as the
  groups of fractions of monoids where
  least common multiples exist, there is no
  nontrivial unit, and some additional finiteness
  conditions are satisfied \cite{Dgk}. Here we
  investigate a wider class of groups of fractions,
  called {\it thin}, which are those associated with
  monoids where minimal common multiples exist,
  but they are not necessarily
  unique. Also, we allow units in the involved
  monoids. The main results are that all thin groups
  of fractions satisfy a quadratic
  isoperimetric inequality, and that, under some
  additional hypotheses, they admit an
  automatic structure.
\end{abstract}

\maketitle

\section{Introduction}

The algebraic theory of braids, as developed in~\cite{Gar}
and~\cite{Eps}, relies on the existence of Garside's
fundamental elements~$\D_n$: for each~$n$, the
braid~$\D_n$ is an element of the monoid~$B_n^+$ which is
a least common multiple of the standard generators~$\s_i$,
and the main technical point is that the left divisors
of~$\D_n$ in~$B_n^+$ coincide with its right divisors. Most
of the results established for Artin's braid groups~$B_n$ have
been extended to more general groups: spherical
Artin groups \cite{Dlg, BrS, Cha}, Garside groups in
the sense of~\cite{Dfx}, and, subsequently,
of~\cite{Dgk} (also called small or thin Gaussian
groups in~\cite{Dfx} and~\cite{Pid}, respectively). All
the considered groups are groups of fractions of
monoids in which least common multiples exist, and,
in each case, a key r\^ole is played by some
element~$\D$ of the associated monoid that satisfies
most of the technical properties of Garside's
braids~$\D_n$. In particular, it is proved
in~\cite{Dgk} that the greedy normal form of braids
\cite{Ady, Eps, ElM} extends to all Garside groups, and
that it gives rise to a bi-automatic structure.

The aim of this paper is to consider groups of fractions of
monoids where common multiples exist, but {\it least}
common multiples need not exist. In this case, no
counterpart of the element~$\D$ need exist in
general, but a number of properties involving the
divisors of~$\D$ can still be established when
considering subsets of the monoid that are closed
under convenient operations. In this
way, one can define an extended notion of normal
form, which coincides with the greedy normal form
when least common multiples exist. The price to pay
for the lack of lcm is a possible non-uniqueness.
However, we shall see that, at least in good cases, this
non necessarily unique normal form is still associated
with an automatic structure. We shall be mainly
interested in the {\it thin} case, defined as the case
when a {\it finite} set of generators with good closure
properties exists. The main results we prove are:

\begin{theo}[Prop.~\ref{P:quis}]\label{T:isop}
  Every thin group of fractions satisfies a quadratic
  isoperimetric inequality.
\end{theo}

\begin{theo}[Prop.~\ref{P:aust}]\label{T:auto}
  Assume that $G$ is the group of fractions of a
  thin cancellative monoid~$M$ that admits a
  Garside element~$\D$ such that all $\D$-normal
  forms in~$M$ have the same length. Then $G$
  is an automatic group.
\end{theo}

These results apply to all thin Gaussian groups,
which are the Garside groups of~\cite{Dgk}, hence in
particular to all spherical Artin groups, for which the
properties were already known, but they also apply to
groups of a completely different flavour, as some
simple examples will show.

The possible interest of our approach is double.
On the one hand, as we mentioned, new groups are
eligible. On the other hand, we hope that extending
classical results may help to understand them better
and to capture the really important hypotheses: 
studying Gaussian groups showed in~\cite{Dgk} that
the fact that Garside's element~$\D_n$ is a least
common multiple of the generators~$\s_i$ of~$B_n$
is useless, and, so, using such a fact gives slightly
misleading arguments. Similarly, the approach
developed in the current paper shows that a clear
distinction should be made between the family of all
divisors of~$\D$ (the ``simple'' elements), and a
smaller subfamily (the ``primitive'' elements) which
contains the real information: the latter can be
extended to the more general framework, while the
former cannot, at least if we use the classical
definition. This leads us here to an alternative,
hopefully improved definition of a simple element. In
the current framework, the proof that the Garside
groups are automatic reduces to a small number of
technical lemmas, each of which is specially easy
when lcm's exist (Lemmas~\ref{L:cov1},  \ref{L:norl},
and~\ref{L:norr}).

The organization of the paper is as follows. In
Sec.~2, we introduce the notion of a spanning subset
of a monoid, which is a generating set satisfying some
additional closure property. Then a thin monoid is
defined to be a monoid that admits a finite spanning
subset. In Sec.~3, we introduce the weaker notion of a
\q-spanning set so as to allow nontrivial units. In
Sec.~4, we define thin groups of fractions as those
associated with a thin Ore monoid, and we prove
Theorem~\ref{T:isop}. In Sec.~5, we show that every
thin monoid admits a minimal spanning subset. In
Sec.~6, we introduce the notion of an $\SS$-simple
element associated with a spanning subset~$\SS$,
which is a counterpart for the notion of divisor
of~$\D_n$ in braid monoids. In Sec.~7, we use
$\SS$-simple elements to construct a counterpart to the
greedy normal form of braids. In Sec.~8, we 
introduce Garside elements, which are convenient
generalizations for the fundamental braids~$\D_n$.
Finally, in Sec.~9, we prove Theorem~\ref{T:auto}.

\section{Spanning subsets of a monoid}\label{S:full}

We consider in the sequel cancellative monoids.
Most of the results until Sec.~4 are valid
if we only assume left cancellativity. If
$M$ is a monoid, we say that an element of~$M$ is
a left (right) unit if it admits a left (right) inverse;
provided $M$ is left or right cancellative, $\uu
\vv = 1$ implies $\uu \vv \uu =\uu$ and $\vv = \vv
\uu \vv$, hence $\vv \uu = 1$, so left and right units
coincide, and they form a subgroup of~$M$ that will
be denoted by~$\UM$. For~$\uu
\in \UM$, we denote by~$\uu\ii$ the (unique) left
and right inverse of~$\uu$. As multiplying by a unit
on the right is often considered in the sequel, we
introduce a notation:

\begin{defi}
  Assume that $M$ is a (left) cancellative monoid. For
  $\xx$, $\xx' \in M$, we say that $\xx \simeq \xx'$
  holds if we have $\xx' = \xx \uu$ for some~$\uu$
  in~$\UM$. We say that a subset~$\SS$ of~$M$ is
  {\it \q-finite} if it contains finitely many
  $\simeq$-classes.
\end{defi}

The relation~$\simeq$ is an equivalence relation
which is compatible with left multiplication. For $\SS
\ince M$, the set~$\SS\op \UM$ is the smallest
$\simeq$-saturated subset of~$M$
including~$\SS$. 

If $M$ is a monoid,  and $\xx$, $\yy$ lie in~$M$, we
say that $\xx$ is a left divisor of~$\yy$, written $\xx \dive
\yy$, if $\yy = \xx \zz$ holds for some~$\zz$. If, in
addition, $\zz$ is not a unit, we say that $\xx$ is a proper
left divisor of~$\yy$, and write $\xx \div \yy$.
We have the symmetric notion of a right divisor,
but, as left divisors play a distinguished r\^ole, we
shall usually simply say ``divisor'' for ``left divisor''.
The set of all (left) divisors of~$\xx$ is
denoted~$\Dl(\xx)$. If $\xx$ is  a (left) divisor
of~$\yy$, we equivalently say that $\yy$ is a (right)
multiple of~$\xx$. Notice that $\simeq$ is
compatible with $\dive$ and $\div$ in the sense that
$\xx \dive \yy$ (\resp $\xx \div \yy$) is equivalent
to $\xx' \dive
\yy'$ (\resp $\xx' \div \yy'$) whenever $\xx'
\simeq \xx$ and $\yy' \simeq \yy$ hold.

The central notion of this paper is that of a spanning
subset of a monoid; it is defined by means of some
closure properties involving left divisors:

\begin{defi}
  Assume that $M$ is a (left) cancellative monoid, and
  $\SS$ is a subset of~$M$. We say that $\SS$ {\it
  spans} $M$ if $\SS$  generates $M$,
  it contains~$1$, it is 
  $\simeq$-saturated, \ie, $\SSU \ince \SS$ holds, 
  and 
  \begin{equation}\label{E:coni}
    \begin{matrix}
      \text{ If we have $\xx \dive \zz$ and $\yy \dive
      \zz$ with $\xx$, $\yy \in \SS$,}\hfill\\
      \text{ then there
       exist $\xx'$, $\yy'$ in~$\SS$ satisfying $\xx\yy' = \yy
       \xx' \dive \zz$.}
     \end{matrix}
   \end{equation}
  We say that $M$ is {\it thin} (\resp {\it \q-thin}) if it
  admits a finite (\resp \q-finite) spanning subset.
\end{defi}

By definition, a thin monoid is finitely generated, but
the converse need not be true, as a spanning subset is more
than a generating subset. Note that the converse of
Implication~\eqref{E:coni} always holds:
$\xx\yy'  = \yy \xx' \dive \zz$ trivially implies $\xx \dive
\zz$ and $\yy \dive \zz$. Spanning subsets always
exist: if $M$ is a monoid, then $M$ is a  spanning
subset of itself. Actually, we shall be mainly interested
in the case when small spanning subsets exist, so, typically,
in the thin case. For a monoid with no nontrivial unit, 
or, more generally, with finitely many units, being
\q-thin is equivalent to being thin.

\begin{exam}\label{X:gaus}
  Let $M$ be a spherical Artin monoid, \ie, one associated
  with a finite Coxeter group~$W$. For $\xx$, $\yy \in M$,
  define $\xx \dr \yy$ to be the unique element~$\zz$
  such that $\xx\zz$ is the least common multiple
  of~$\xx$ and~$\yy$. Then the closure of the standard
  generators~$\s_i$ under the operation~$\dr$  is a finite
  spanning subset of~$M$, in one-to-one correspondence
  with some subset of~$W$ \cite{Dfx}. So spherical
  Artin monoids are thin.
  
  More generally, if $M$ is a Gaussian monoid in the
  sense of~\cite{Dfx, Dgk}, \ie, a cancellative monoid 
  in which least common multiples exist and division
  has no infinite descending chain, then the closure of
  the set of atoms under operation~$\dr$ is a
  minimal spanning subset of~$M$. Thus the
  monoid~$M$ is thin if and only if the latter closure
  is finite, \ie, if $M$ is thin in the sense of~\cite{Pid}
  (or small in the sense of~\cite{Dfx}), so the
  terminologies are compatible.\footnote{In order
  to uniformize terminology with other authors,
  we use {\it Garside monoid} as a synonym for
  thin Gaussian monoid, and {\it Garside group} for
  the group of fractions of a Garside
  monoid~\cite{Dgk}.}
\end{exam}

\begin{figure}
  \includegraphics{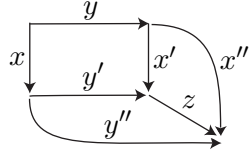}
  \caption{Spanning subset}\label{F:full}
\end{figure}

\begin{lemm}
  Assume that $M$ is a (left) cancellative monoid, and
  $\SS$ is subset of~$M$. Then 
  Condition~\eqref{E:coni} for~$\SS$
  is equivalent to
  \begin{equation}\label{E:conj}
    \begin{matrix}
      \text{If we have $\xx\yy'' = \yy\xx''$ with $\xx$, $\yy \in
      \SS$, then there exist $\xx'$, $\yy'$
      in~$\SS$}\hfill\\
      \text{and $\zz$ in~$M$ satisfying $\xx\yy' = \yy\xx'$, $\xx''
      = \xx' \zz$, and $\yy'' = \yy' \zz$
       (Fig.~\ref{F:full}).}
     \end{matrix}
   \end{equation}
\end{lemm}

\begin{proof}
   It is clear that $\eqref{E:conj}$ implies~\eqref{E:coni};
conversely,
   assuming $\xx\yy'' = \yy\xx''$, Condition~\eqref{E:coni}
   implies that there exist $\xx'$, $\yy'$ in~$\SS$ and
   $\zz$ in~$M$ satisfying $\xx\yy' = \yy\xx'$, $\xx\yy''=
   \xx\yy'\zz$, and $\yy\xx'' = \yy\xx'\zz$, hence
   $\xx'' = \xx' \zz$ and $\yy'' = \yy' \zz$ if we
   can cancel $\xx$ and $\yy$ on the left.
\end{proof}

If $\SS$, $\TT$ are subsets of a monoid~$M$, we put
$\SS \op \TT = \{\xx\yy \, ; \, \xx \in \SS, \yy \in \TT\}$.
In particular, $\SS^p$ is the set of all elements
that can be written as $\xx_1 \cdots \xx_p$ with
$\xx_1$, \pp,
$\xx_p \in \SS$. Notice that  $1 \in \SS$
implies $\SS \ince \SS^p$ for~$p \ge 2$. We put
$\SS^0 = \{1\}$. 

\begin{lemm}\label{L:iter}
  Assume that $M$ is a monoid, and $\SS$ is a subset
  of~$M$ satisfying Condition~\eqref{E:conj}. Then, if we have 
  $\xx\yy'' = \yy\xx''$ with $\xx \in \SS^p$ and
  $\yy \in \SS^q$, there exist an element~$\zz$ of~$M$
  and two sequences $\xx_{i, j}$, $\yy_{i, j}$, $0 \le i \le p$,
  $0 \le j \le q$, of elements of~$\SS$ satisfying
  $\xx_{i, j-1}\yy_{i,j} = \yy_{i-1, j}\xx_{i, j}$ for
  all~$i$, $j$, and $\xx = \prod \xx_{i,0}$, $\yy = \prod
  \yy_{0,j}$, $\xx'' =  \prod \xx_{i, q} \, \zz$, and $\yy'' =
  \prod \yy_{p, j} \, \zz$. So, in particular, there exist
  $\xx'$ in~$\SS^p$, $\yy'$ in~$\SS^q$ and $\zz$
  in~$M$ satisfying $\xx \yy' = \yy \xx'$, $\xx'' = \xx'
  \zz$, and $\yy'' = \yy' \zz$.
\end{lemm}

\begin{figure}
  \includegraphics{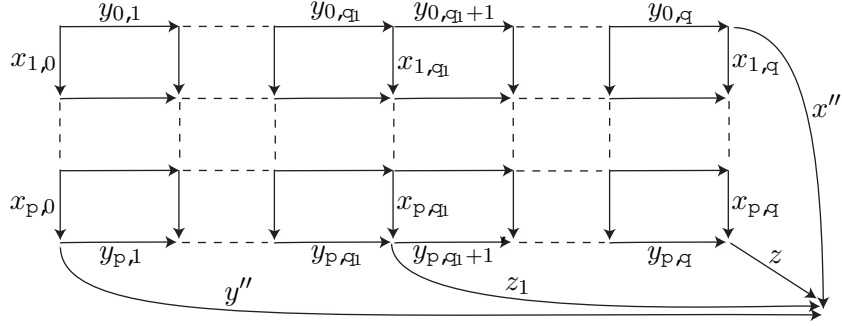}
  \caption{Power of a spanning subset}\label{F:powe}
\end{figure}

\begin{proof}
  (Fig.~\ref{F:powe})
  First, the condition is sufficient, as the local equalities
  $\xx_{i,j} \yy_{i,j+1} = \yy_{i, j} \xx_{i+1, j}$ imply
  $\prod_i \xx_{i, 0} \prod_j \yy_{p, j} = \prod_j \yy_{0, j}
  \prod_i \xx_{i, q}$, hence $\xx \yy'' = \yy \xx''$ when
  $\xx$, $\yy$, $\xx''$, $\yy''$ have the above specified
  values.
  
  We prove now that the condition is necessary. The result is
  trivial for $p = 0$ or $q = 0$. Indeed, $p = 0$ means
  $\xx = 1$:
  then the hypothesis $\yy \in \SS^q$ allows us to write
  $\yy = \prod_j \yy_{0,j}$ with $\yy_{0, 1}$, \pp, $\yy_{0, q}
  \in \SS$, and the hypothesis $\yy'' = \yy \xx''$ then
  gives
  $\yy'' = \prod_j \yy_{0,j} \zz$ with $\zz = \xx''$.
  Then we use induction on~$p + q$.  By the remark
  above, the
  first nontrivial case is $p = q = 1$, and, then, the result is
  true by Condition~\eqref{E:conj}.
  Assume now $p + q  \ge 3$, with $p, q \ge 1$. Then at least
  one of $p$, $q$ is greater than~$1$.  Assume for instance
  $q \ge 2$. Write
  $\yy = \yy_1 \yy_2$ with $\yy_1 \in \SS^{q_1}$, $\yy_2 \in
  \SS^{q_2}$ and $1 \le q_1, q_2 < q$. 
  Applying the induction hypothesis to $\xx \in \SS^p$, $\yy_1
  \in \SS^{q_1}$ and $\xx \yy'' = \yy_1 (\yy_2 \xx'')$
  gives
  $\zz_1$ in~$M$ and $\xx_{i, j}$, $\yy_{i, j}$, $0 \le i
  \le p$, $0 \le j \le q_1$ in~$\SS$ satisfying
  $\xx_{i, j-1}\yy_{i,j} = \yy_{i-1, j}\xx_{i, j}$, $\yy_2
  \xx'' = \xx_1 \zz_1$, $\yy'' = \yy'_1
  \zz_1$ with 
  \[
    \xx = \prod_1^p \xx_{i,0}, \quad
    \yy_1 = \prod_1^{q_1} \yy_{0,j}, \quad
    \xx_1 =  \prod_1^p \xx_{i, q_1} , \quad
    \yy'_1 = \prod_1^{q_1} \yy_{p, j}.
  \]
  Applying the induction hypothesis to~$\xx_1 \in \SS^p$,
  $\yy_2 \in \SS^{q_2}$ and $\xx_1 \zz_1 = \yy_2 \xx''$
  gives
  $\zz$ in~$M$ and $\xx_{i, j}$, $\yy_{i, j}$, $0 \le i
  \le p$, $ q_1 <  q_2$ in~$\SS$ satisfying
  $\xx_{i, j-1}\yy_{i,j} = \yy_{i-1, j}\xx_{i, j}$, $ \xx''
  = \xx' \zz$, $\zz_1 = \yy'_2 \zz$ with 
  \[
    \xx = \prod_1^p \xx_{i,0}, \quad
    \yy_2 = \prod_{q_1+1}^{q_2} \yy_{0,j}, \quad
    \xx' =  \prod_1^p \xx_{i, q} , \quad
    \yy'_2 = \prod_{q_1+1}^{q_2} \yy_{p, j}.
  \]
  Putting $\yy' = \yy'_1 \yy'_2$ gives the expected result.

  Finally, with the previous notation, put $\xx' =\prod_i 
  \xx_{i, q}$ and $\yy' = \prod_j \yy_{p, j}$. By 
  construction, $\xx'$ lies in~$\SS^p$, $\yy'$ lies
  in~$\SS^q$, and we have $\xx \yy' = \yy \xx'$, $\xx'' =
  \xx' \zz$, and $\yy'' = \yy' \zz$.
\end{proof}

Applying Lemma~\ref{L:iter} with $q = p$, we obtain

\begin{prop}\label{P:powe}
 Assume that $M$ is a (left) cancellative monoid. If
  $\SS$ spans~$M$, so does $\SS^p$ for
 every positive~$p$.
\end{prop}

Cancellativity is not used in the proof of
Lemma~\ref{L:iter}, so, at the expense of using \eqref{E:conj} instead
of~\eqref{E:coni} in the definition of a spanning subset, we could state
Prop.~\ref{P:powe} for a general monoid.

\begin{prop}\label{P:clos}
  Assume that $M$ is a (left) cancellative monoid, and
  $\SS$ spans $M$. Then every right divisor
  of an element of~$\SS$ lies in~$\SS$, and we
  have $\UM \op \SS \ince \SS$.
\end{prop}

\begin{proof}
  Assume $\yy = \xx\zz \in \SS$. As $\SS$
  generates~$M$, we have $\xx \in \SS^p$ for
  some~$p$, so, applying
  Lemma~\ref{L:iter} to the equality $\xx \cdot \zz = 
  \yy \cdot 1$, we find $\xx$ in~$\SS^p$, $\yy'$
  in~$\SS$, and $\zz'$ in~$M$ satisfying $\xx\yy' =
  \yy\xx'$, $\zz = \yy'\zz'$, and $1 = \xx'\zz'$. The
  latter relation shows that $\xx'$ and $\zz'$ are units
  and $\zz = \yy'\zz'$ then implies $\zz \in \SS \op
  \UM$, hence $\zz \in \SS$ as $\SS$ is supposed
  to be closed under right multiplication by a unit.
  
  Assume $\xx \in \SS$ and $\uu \in \UM$. Then we have
  $\xx = \uu\ii (\uu \xx)$, so $\uu \xx$ is a right divisor
  of~$\xx$, and, by the previous result, it belongs
  to~$\SS$.
\end{proof}

One of the interests of spanning subsets is that they
completely determine the monoid in the sense
below. In the sequel, if $\SS$ is a set (of letters), and
$R$ is a set of relations over~$\SS$, \ie, of equalities
of the form $\xx_1 \cdots \xx_p = \yy_1 \cdots
\yy_q$ with
$\xx_1$, \pp, $\yy_q \in \SS$, we denote by
$\Mon(\SS; R)$ the monoid so presented, and by
$\Gr(\SS; R)$ the group with the same presentation.

\begin{defi}
  Assume that $M$ is a (left) cancellative monoid, and
  $\SS$ spans~$M$. We denote by $\RX$ the set
  of all relations $\xx\yy' = \yy\xx'$ with $\xx$, $\yy$,
  $\xx'$, $\yy' \in \SS$. 
\end{defi}

\begin{prop}\label{P:pres}
  Assume that $M$ is a (left) cancellative monoid, and
  $\SS$ spans $M$. Then $\Mon(\SS; \RX)$ is a
  presentation
  of~$M$, and every equality $\xx_1 \cdots \xx_p = 
  \yy_1 \cdots \yy_q$ with $\xx_1$, \pp, $\yy_q \in
  \SS$ can be proved using $O((p+q)^2)$ relations
  of~$\RX$.
\end{prop}

\begin{figure}
  \includegraphics{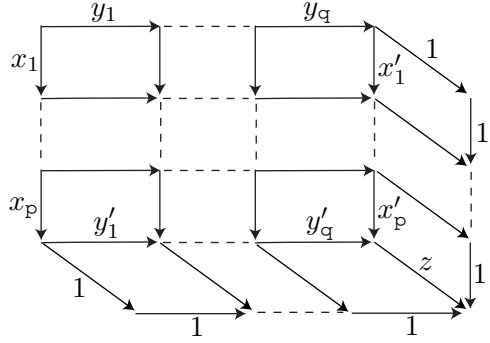}
   \caption{Quadratic isoperimetric inequality}
  \label{F:pres}
\end{figure}

\begin{proof} (Fig. \ref{F:pres})
  The set~$\SS$ generates~$M$ by definition. Assume
  $\xx_1 \cdots \xx_p = \yy_1
  \cdots \yy_q$ with $\xx_1$, \pp, $\yy_q \in \SS$.
By Lemma~\ref{L:iter}, there
  exist $\xx'_1$, \pp, $\xx'_p$, $\yy'_1$, \pp, $\yy'_q
\in \SS$ and $\zz \in M$ satisfying
  \[
    \xx_1 \cdots \xx_p \yy'_1 \cdots \yy'_q = \yy_1 \cdots \yy_q
    \xx'_1 \cdots \xx'_p, \quad
    \xx'_1 \cdots \xx'_p \zz = 1,  \quad\text{and }
    \yy'_1 \cdots \yy'_q \zz = 1,
  \]
  and, moreover, the first equality can be established using
  $pq$~relations in~$\RX$. As for the other ones, we
  know by Prop.~\ref{P:clos} that each of the elements
  $\xx'_i \cdots \xx'_p\zz$ and $\yy'_j \cdots \yy'_q \zz$
  belongs to~$\SS$, and, therefore, the equality
  $\xx'_1 \cdots \xx'_p \zz = 1$ can be established using
  $p$~relations of~$\RX$ (of the special form $\xx\yy =
  1$), and, similarly, $\yy'_1 \cdots \yy'_q \zz = 1$ can be
  established using $q$~relations of~$\RX$. So, finally, the
  equality $\xx_1 \cdots \xx_p = \yy_1
  \cdots \yy_q$ can be established using at most $(p+q)^2/4
  + (p+q)$ relations of~$\RX$.
\end{proof}

Applying the previous result to the case of a finite
spanning subset, we obtain:

\begin{prop}
  Every thin cancellative monoid satisfies a quadratic
  iso\-perimetric inequality.
\end{prop}
  
In order to construct new examples of thin monoids,
Prop.~\ref{P:pres} suggests that we consider
presentations where the relations are of the form
$\xx\yy' =
\yy\xx'$, \ie, involve words of length~$2$ at most.
Every monoid admits a presentation of this type,  and
the question arises of recognizing spanning subsets.
Typically, if $\Mon(\SS; R)$ is a presentation of the
type above for a monoid~$M$, there is no obvious
reason  why $\SS$ should span~$M$, as some
equalities in~$M$ may follow from the relations
of~$R$ but not decompose into such relations using
the scheme of Fig.~\ref{F:pres}. In particular, there is
no reason  why the equality $R = \RX$ should hold.
Here, we shall refer to~\cite{Dgp}, where the
notion of a complete presentation is defined. The idea
is that a presentation is complete if the relations have
no hidden consequence,
\ie, if enough relations have been displayed to
avoid any such
hidden consequence. Then the result is that, if
$\Mon(\SS; R)$ is a complete presentation for the
monoid~$M$, then a sufficient condition for $\SS$ to 
span~$M$ is that each relation in~$R$ has
length~$2$ at most.

\begin{exam}\label{X:main}
  With this method, we can exhibit thin monoids that
  do not resemble the Gaussian monoids of
  Example~\ref{X:gaus}, namely
  monoids where least common multiples do not exist. The
  following three examples are typical, and they will be
  considered throughout the paper:
  \begin{gather*}
    M_1 = \Mon(\PresI), \\
    M_2 = \Mon(\PresII), \\
    M_3 = \Mon(\PresIII).
  \end{gather*}
  Applying the criterion of~\cite{Dgp}, one checks
  that the above presentations are complete, and that, in
  each case, the involved set of generators completed
  with~$1$ is a spanning subset. Thus the monoids~$M_i$
  are thin. Moreover, there is no relation $\xx\yy =
  \xx\yy'$ or $\yy\xx = \yy'\xx$
  with $\yy \not= \yy'$ in the above presentations, so,
  according to~\cite{Dgp} again, the monoids~$M_i$ are
  cancellative.
\end{exam}

\begin{figure}
  \includegraphics{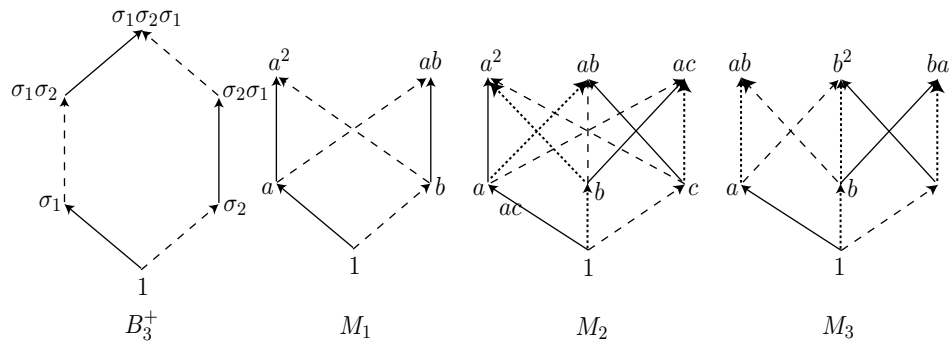}
  \caption{Characteristic graph associated 
  with a spanning subset}\label{F:grap}
\end{figure}

If $\SS$ spans a monoid~$M$, then
$M$ is determined by~$\SS$
and~$\RX$. Especially when the considered
set~$\SS$ is finite, \ie, in the thin case, it is
natural to introduce the subgraph of the Cayley
graph of~$M$ displaying the relations of~$\RX$:
by the remarks above, such a (finite) graph
completely determines the monoid. In the
Gaussian case, \ie, when least common multiples
exist, the graph is a lattice, in the sense that any
two vertices admit a unique immediate common
successor. In the general case, this need not be
true. For instance, we display in Fig.~\ref{F:grap}
the graphs associated with the braid
monoid~$B_3^+$ and the spanning subset $\{1,
\s_1, \s_2$, $\s_1\s_2$, $\s_2\s_1\}$, and with the
thin monoids~$M_1$, $M_2$, $M_3$ of
Example~\ref{X:main}.

\section{\Q-spanning subsets}

Many results for the thin case extend to the
\q-thin case, and, to deal with the latter, it is
convenient to introduce the notion of a \q-spanning subset
of a monoid, of which a typical example is a
$\simeq$-selector through a spanning set, \ie, a
subset that picks one element in each equivalence
class.

\begin{defi}
  Assume that $M$ is a (left) cancellative monoid,
  and $\SS$ is a subset of~$M$. We say that $\SS$ {\it
  \q-spans}~$M$ if $\SS \op \UM$ spans~$M$. 
\end{defi}

By definition, a spanning subset is \q-spanning, and both
notions coincide if there is no nontrivial unit.

\begin{prop}\label{P:sele}
  Assume that $M$ is a (left) cancellative monoid,
  and $\SS$ spans $M$. Then every
  $\simeq$-selector through~$\SS$ \q-spans~$M$. Conversely, if $\SS$ is a minimal
  \q-spanning subset of~$M$, then $\SS$ is a
  $\simeq$-selector.
\end{prop}

\begin{proof}
  If $\S$ is a $\simeq$-selector through~$\SS$,
  we have $\S \op \UM = \SS$, so $\S$
  is \q-spanning. 
  On the other hand, assume that $\SS$ is a
  minimal \q-spanning subset of~$M$.
  Assume $\yy \simeq \xx$, with $\xx$, $\yy \in
  \SS$ and $\xx \not= \yy$.  Then 
  $(\SS - \{\yy\}) \op \UM$ is equal to~$\SSU$,
  and, therefore, $\SS - \{\yy\}$ \q-spans~$M$, 
  contradicting the minimality of~$\SS$.
\end{proof}

\begin{coro}
  A cancellative monoid is \q-thin if and only if it
  admits a finite \q-spanning subset.
\end{coro}

\begin{lemm}\label{L:quas}
  Assume that $M$ is a (left) cancellative monoid,
  and $\SS$ is a subset of~$M$. Then $\SS$
  \q-spans~$M$ if and only if $\SS$
  contains~$1$, $\SS\op\UM$ generates $M$, $\UM \op
  \SS \ince \SS \op \UM$ holds, and
  \begin{equation}\label{E:conk}
    \begin{matrix}
      \text{If we have $\xx \dive \zz$ and $\yy \dive
      \zz$ with $\xx$, $\yy \in \SS$,}\hfill\\
      \text{then there
       exist $\xx'$, $\yy'$ in~$\SS$ satisfying $\xx\yy'
       \simeq \yy \xx' \dive \zz$.}
     \end{matrix}
   \end{equation}
\end{lemm}

\begin{proof}
  Assume that $\SS$ satisfies the above conditions.
  Then $\SSU$ contains~$1$, it generates~$M$ by
  hypothesis,
  and it is $\simeq$-saturated by construction. Assume
  $\xx \uu \dive \zz$ and $\yy \vv \dive \zz$ with $\xx$,
  $\yy \in \SS$, and $\uu$, $\vv \in \UM$.
  Then we have $\xx \dive \zz$ and $\yy \dive \zz$,
  so, by~\eqref{E:conk}, there exist $\xx'$, $\yy'$ in~$\SS$
  satisfying $\xx\yy' \simeq \yy \xx' \dive \zz$, say
  $\xx\yy' = \yy\xx'\ww \dive \zz$ with $\ww \in \UM$.
  Then we have also $(\xx\uu) (\uu\ii\yy') =
  (\yy\vv) (\vv\ii \xx'\ww) \dive \zz$, and the elements
  $\uu\ii\yy'$ and $\vv\ii \xx'\ww$ belong to $\UM
  \op \SS \op \UM$, hence to~$\SSU$, which
  therefore
  satisfies Condition~\eqref{E:coni}, and spans~$M$.
  
  Conversely, assume that $\SS$ \q-spans~$M$. Then $1 \in \UM$ implies $\SS
  \ince \SSU$, and restricting Condition~\eqref{E:coni}
  for~$\SSU$ to~$\SS$ yields Condition~\eqref{E:conk}.
  Moreover, we have $\UM \op \SSU \ince \SSU$
  by Prop.~\ref{P:clos}, hence, a fortiori,
  $\UM \op \SS \ince \SSU$. So all conditions of 
  Lemma~\ref{L:quas} are satisfied.
\end{proof}

If $\SS$ \q-spans a monoid~$M$, then 
we have $\UM \op \SS \ince \SSU$ by 
Lemma~\ref{L:quas}, and
a straightforward induction then implies
\begin{equation}\label{E:poqu}
  (\SSU)^n = \SS^n \op \UM
\end{equation}
for every positive~$n$. As $\SSU$ is supposed to
generate~$M$, it follows that every element of~$M$
admits a decomposition of the form $\xx = \xx_1
\cdots \xx_n \uu$ with $\xx_1$, \pp, $\xx_n \in
\SS$ and $\uu \in \UM$.

\begin{prop}\label{P:powf}
  Assume that $M$ is a (left) cancellative monoid and 
  $\SS$ \q-spans~$M$. 
  
  \noindent (i) The subset~$\SS^p$ \q-spans~$M$
  for every positive~$p$.

  \noindent (ii) If we have $\xx
  \dive \zz$ and $\yy \dive \zz$ with $\xx \in \SS^p$
  and $\yy \in \SS^q$, then there exist $\xx' \in \SS^p$
  and $\yy' \in \SS^q$ satisfying $\xx\yy' \simeq \yy\xx'
  \dive \zz$.
  
  \noindent (iii) If we have $\xx_i
  \dive \zz$ with $\xx_i \in \SS$ for $1 \le i \le n$, then
  there exists $\xx$ in~$\SS^n$ satisfying $\xx_i
  \dive \xx \dive \zz$ for $1 \le i \le n$. 
\end{prop}

\begin{proof}
  (i) By Lemma~\ref{L:quas}, $\SSU$ spans~$M$, so,
by Prop.~\ref{P:powe}, $(\SSU)^p$ spans~$M$ as well
for every positive~$p$.
  By~\eqref{E:poqu}, the latter is~$\SS^p \op \UM$,
  and, by Lemma~\ref{L:quas} again, $\SS^p
  \op \UM$ spanning~$M$ implies $\SS^p$ 
  \q-spanning~$M$.

  (ii) Applying Lemma~\ref{L:iter} to the spanning
  subset~$\SSU$ of~$M$,
  we obtain $\xx''$ in $(\SSU)^p$ and $\yy''$
  in~$(\SSU)^q$ satisfying $\xx\yy'' = \yy\xx''
  \dive \zz$. By~\eqref{E:poqu}, we have $(\SSU)^p =
  \SS^p \op \UM$ and $(\SSU)^q = \SS^q \op \UM$, so we
  deduce that there exist $\xx'$ in~$\SS^p$, $\yy'$
  in~$\SS^q$, and $\uu$, $\vv$ in~$\UM$
  satisfying
  $\xx'' = \xx'\uu$ and $\yy'' = \yy' \vv$, hence $\xx\yy'
  \simeq \yy\xx' \dive \zz$.

  (iii) We use induction on~$n \ge 1$. For $n = 1$, we
  take $\xx = \xx_1$.  Assume $n \ge 2$. By induction
  hypothesis, some $\yy$ in $\SS^{n-1}$
  satisfies $\xx_i \dive \yy \dive \zz$ for $i \le n-1$.
  Applying~(ii) to~$\xx_n$ and~$\yy$, 
  we obtain $\xx'$ in~$\SS$ and $\yy'$ in~$\SS^{n-1}$
  satisfying $\xx_n \yy' \simeq \yy \xx' \dive \zz$.
  Putting $\xx = \yy \xx'$ gives the result.
\end{proof}

\begin{prop}
  Assume that $M$ is a (left) cancellative monoid and $\SS$ is a
  minimal \q-spanning subset of~$M$. Then, for every
  pair~$(\xx, \uu)$ in~$\SS \times
  \UM$, there exists a unique pair~$(\xx', \uu')$
  in~$\SS \times \UM$ satisfying $\uu\xx =
  \xx'\uu'$. The mapping $\xx \mapsto \xx'$
  defines an action of the group~$\UM$ on~$\SS$. If
  $\xx$ is fixed under this action, the mapping $\uu
  \mapsto \uu'$ is an endomorphism of~$\UM$ for
  every~$\xx$.
\end{prop}

\begin{proof}
  By Lemma~\ref{L:quas}, we have $\UM \op \SS \ince
  \SSU$, so $\uu \xx \in \SSU$, \ie, there
  exist $\xx'$ in~$\SS$ and $\uu'$ in~$\UM$
  satisfying $\uu \xx =\xx' \uu'$. The  uniqueness
  of~$\xx'$ follows from~Prop.~\ref{P:sele}, and that
  of~$\uu'$ then follows from
  left cancellativity. Writing ${}^\uu\xx$ for~$\xx'$,
  we have $1\xx
  \simeq \xx$ for every~$\xx$ in~$\SS$, and, therefore
  ${}^1\xx = \xx$, and $\uu\vv\xx \simeq \uu {}^\vv\xx
  \simeq {}^\uu({}^\vv\xx)$, hence ${}^{\uu\vv}\xx =
  {}^\uu({}^\vv\xx)$.
  
  Assume now ${}^\uu \xx = \xx$. Writing $\uu^\xx$
  for~$\uu'$, we find $(\uu\vv)\xx = \uu \xx
  \vv^\xx = \xx \uu^\xx \vv^\xx$, hence $(\uu
  \vv)^\xx = \uu^\xx \vv^\xx$.
\end{proof}

The previous result gives a way for constructing
\q-thin monoids with a prescribed
group of units. Assume that $M$ is a cancellative monoid
with no nontrivial unit, $\SS$ spans~$M$, and
$G$ is a group with a left action on~$M$ that preserves
$\SS$ globally. Consider the semi-direct product $M
\rtimes G$ where $(\xx, \uu)(\yy, \vv)$ is defined to be
$(\xx\uu(\yy), \uu\vv)$. The set of units in $M
\rtimes G$ is $\{1\} \times G$, and the set $\SS \times
\{1\}$ is a \q-spanning subset of $M \rtimes G$. More
generally, instead of a semidirect product, we
could also use a crossed product as in~\cite{Pid}.

In the same way as a spanning subset determines a
monoid, a \q-spanning subset together with the units
determine a monoid.

\begin{defi}
  Assume that $M$ is a (left) cancellative monoid, and
  $\SS$ \q-spans~$M$. We denote by
  $\RRX$ the set consisting of
  
  (i) all relations $\xx\yy' = \yy\xx'\uu$ with $\xx$,
  $\yy$, $\xx'$, $\yy' \in \SS$ and $\uu \in \UM$,
  
  (ii) all relations $\uu \xx = \xx' \uu'$ with $\xx$, $\xx'
  \in \SS$ and $\uu$, $\uu' \in \UM$,
  
  (iii) all relations $\uu \vv = \ww$ with $\uu$, $\vv$, $\ww
  \in \UM$.
\end{defi}

\begin{prop}\label{P:pret}
  Assume that $M$ is a (left) cancellative monoid, and
  $\SS$ \q-spans~$M$. Then $\Mon(\SS
  \cup \UM; \RRX)$ is a
  presentation of~$M$, and every equality $\xx_1 \cdots \xx_p = 
  \yy_1 \cdots \yy_q$ with $\xx_1$, \pp, $\yy_q \in
  \SS$ can be
  proved using $O((p+q)^2)$ relations of~$\RRX$.
\end{prop}

\begin{figure}
  \includegraphics{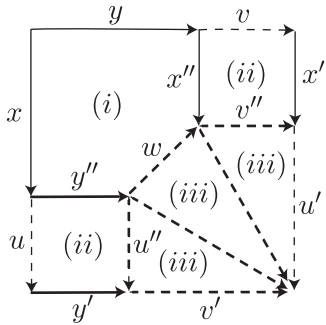}
  \caption{Decomposition of a relation}\label{F:pret}
\end{figure}

\begin{proof}
  By Lemma~\ref{L:quas}, the set~$\SSU$ spans~$M$,
  so, by Prop.~\ref{P:pres},  $\Mon(\SSU; R_{\SS \op
  \UM})$ is a
  presentation of~$M$, and every equality $\xx_1 \cdots
  \xx_p = \yy_1 \cdots \yy_q$ with $\xx_1$, \pp, $\yy_q
  \in \SSU$ can be established using at most
  $O((p+q)^2)$ relations in~$R_{\SS \op \UM}$. So, it
  suffices to prove that every relation
  in~$R_{\SS \op \UM}$ can be decomposed into a
  uniformly bounded number of relations in~$\RRX$. By
  construction,
  every element in~$\SSU$ can be expressed as $\xx
  \uu$ with $\xx \in \SS$ and $\uu \in \UM$. 
  So assume $\xx \uu \yy' \vv'= \yy \vv \xx' \uu'$ with
  $\xx$, $\yy$, $\xx'$, $\yy' \in \SS$ and $\uu$,
  $\vv$, $\uu'$, $\vv'\in \UM$. There exist
  $\xx''$, $\yy''$ in~$\SS$ and $\uu''$, $\vv''$
in~$\UM$
  satisfying the type~(ii) relations $\uu \yy' = \yy'' \uu''$
  and $\vv \xx' = \xx'' \vv''$: then $\xx \uu \yy' \vv'= \yy
  \vv \xx' \uu'$  implies the type~(i) relation $\xx \yy'' \ww
  = \yy \xx''$ with $\ww  = \uu'' \vv' {\uu'}\ii {\vv''}\ii$,
  and the latter equality follows from three type~(iii)
  relations (Fig.~\ref{F:pret}). Thus every relation
  of~$R_{\SS \op \UM}$ can be decomposed into at most
  six relations in~$\RRX$.
\end{proof}

Instead of using all units in the presentation, we can
replace $\UM$ with any set that generates it (as a
monoid). We obtain:

\begin{coro}\label{P:presIII}
  Assume that $M$ is a (left) cancellative monoid, 
  $\SS$ \q-spans~$M$, and $\Mon(\SS'; R')$
  is a presentation of~$\UM$ as a monoid.
  Let $R$ consist of
  
  \noindent (i) all relations $\xx\yy' = \yy\xx'\uu_1 
  \cdots \uu_p$ with $\xx$, $\yy$, $\xx'$, $\yy' \in
  \SS$, $\uu_1$, \pp, $\uu_p \in \SS'$,
  
  \noindent  (ii) all relations $\uu \xx = \xx' \uu'_1 
  \cdots \uu'_p$ with $\xx$, $\xx' \in \SS$, $\uu$,
  $\uu'_1$, \pp, $\uu'_p \in \SS'$.
  
  \noindent Then $\Mon(\SS \cup \SS'; R \cup R')$ is a
  presentation of~$M$. 
\end{coro}

The question of whether an isoperimetric inequality is
satisfied when the previous presentation is finite is
open: even if we assume that the presentation of the 
group~$\UM$ satisfies such a condition, there is no 
easy way to conclude for~$M$ as we know nothing
about the elements denoted $\uu''$,
$\vv''$, and~$\ww$ in Fig.~\ref{F:pret}.

\section{Thin groups of fractions}

By a well known result of Ore \cite{CPr}, a cancellative
monoid~$M$ embeds in a group of (right) fractions if and
only if any two elements of~$M$ admit at least one
common multiple. A monoid satisfying such conditions
will be called a {\it Ore monoid} in the sequel. It is
then natural to consider those groups of fractions that
are associated with thin monoids:

\begin{defi}
  We say that a group~$G$ is a {\it thin} (\resp {\it
  \q-thin}) group of fractions if there exists a thin
  (\resp \q-thin) Ore monoid~$M$ such that $G$ is
  the group of fractions of~$M$.
\end{defi}

Thus, the braid groups~$B_n$, the spherical Artin groups,
and, more generally, the Garside groups
of~\cite{Dgk} are thin groups of fractions. We shall give
more examples below.

A nice point is that, when a \q-spanning subset is known in
a monoid~$M$, then it is easy to study the possible
existence of common multiples in~$M$. Indeed, we have
the following criterion:

\begin{prop}\label{P:excm}
  Assume that $M$ is a (left) cancellative monoid, and
  $\SS$ is a \q-spanning subset in~$M$. Then any
  two elements of~$M$ admit a
  common multiple if and only the following
  condition holds:
  \begin{equation}\label{E:excm}
    \text{For all $\xx$, $\yy$ in~$\SS$, there exist
    $\xx'$, $\yy'$ in~$\SS$ satisfying $\xx\yy'
    \simeq \yy\xx'$.}
  \end{equation}
\end{prop}

\begin{proof}
  Assume $\xx$, $\yy \in \SS$. If common
  multiples always exist in~$M$, there 
  exists~$\zz$ satisfying $\xx \dive \zz$ and $\yy \dive
  \zz$, and, therefore, since $\SS$ is \q-spanning, there
  exist $\xx'$, $\yy'$ in~$\SS$ satisfying $\xx\yy' \simeq
  \yy\xx' (\dive \zz)$. So \eqref{E:excm} holds.
  
  Conversely, assume~\eqref{E:excm}. First, we
claim that, for all
  $\xx$, $\yy \in \SSU$, there exist $\xx'$, $\yy'
  \in \SSU$ satisfying $\xx\yy' = \yy\xx'$. Indeed,
  assume $\xx \simeq \xx_0$, $\yy \simeq \yy_0$
  with $\xx$, $\yy \in \SS$. By hypothesis, there
  exists~$\zz$ satisfying $\xx_0 \dive \zz$ and $\yy_0
  \dive \zz$. Then $\xx \dive \zz$ and $\yy \dive \zz$
  hold as well. As $\SSU$ spans~$M$, we
  deduce that there exist $\xx'$, $\yy'$ in~$\SSU$
  satisfying $\xx\yy' = \yy\xx' (\dive \zz)$.
  
  We prove now that, if $\xx$ belongs
  to~$(\SSU)^p$ and $\yy$ belongs to~$(\SSU)^q$, then
  $\xx$ and $\yy$ admit a common multiple
  in~$M$. The result is trivial for $p = 0$ and $q = 0$, and
  the case $p = q = 1$ has been treated above.  Then we use
  a recurrence on~$p + q$: the principle is to
  construct a diagram like the one in Fig.~\ref{F:powe}
  starting from the left and the top edges. By hypothesis,
  each small square can be closed, so, inductively, the 
  full diagram can be completed.
\end{proof}

It follows that, if $M$ is a cancellative monoid and $\SS$ is
a \q-spanning subset in~$M$ that
satisfies~\eqref{E:excm}, then $M$ is a Ore
monoid, and it embeds in a group of right
fractions. Moreover, the latter admits the
presentation~$\Gr(\SS;
\RX)$, where $\RX$ is as in Prop.~\ref{P:pres}. In
particular, we can state: 

\begin{coro}
  Assume that $M$ is a thin cancellative monoid and $\SS$
  spans~$M$ and satisfies~\eqref{E:excm}.
  Then the group~$\Gr(\SS; \RX)$ is a thin group
  of fractions.
\end{coro}

\begin{exam}\label{X:maio}
  For $i = 1, \pp, 3$, the monoid~$M_i$ of
  Example~\ref{X:main} satisfies
  Condition~\eqref{E:excm}, as we can check
  on the graph of Fig.~\ref{F:grap}. Thus, any two elements
  in~$M_i$ admit a common multiple, and,
  therefore, $M_i$ embeds in a group of fractions.
  So the groups
  \begin{gather*}
    G_1 = \Gr(\PresI) \\
    G_2 = \Gr(\PresII) \\
    G_3 = \Gr(\PresIII)
  \end{gather*}
  are thin groups of fractions.
  
  Observe that the monoid $\Mon(\a, \b, \c \, ; \,
  \a\b = \b\c = \c\a)$ is the Birman-Ko-Lee monoid
  for the braid group~$B_3$ \cite{BKL}; then $M_2$ is
  the quotient of the latter
  monoid under the additional relation $\a^2 =
  \b^2$, and,  therefore, the group~$G_2$ is the
  quotient of~$B_3$ obtained by
  adding the relation $\s_1^2 = \s_2^2$, thus an
  intermediate group between~$B_3$ and the
  symmetric group~$S_3$.
\end{exam}

We can now state our first general result about thin groups
of fractions:

\begin{prop}\label{P:quis}
  Every thin group of fractions satisfies a quadratic
  isoperimetric inequality.
\end{prop}

\begin{proof}
  Assume that $G$ is the group of fractions of the Ore
  monoid~$M$, and $\SS$ is a finite spanning subset
  of~$M$. By Prop.~\ref{P:pres}, $\Mon(\SS; \RX)$
  is a presentation of~$M$, and, therefore, 
  $\Gr(\SS; \RX)$ is a presentation of~$G$, which is finite
  by construction.
  Assume that $\xx_1^{e_1} \cdots \xx_n^{e_n} = 1$
  holds in~$G$, with $\xx_1$, \pp, $\xx_n \in \SS$
  and $e_1, \pp, e_n = \pm 1$. First, an easy
  induction on~$m$ shows that, if the sequence $(e_1,
  \pp, e_n)$ contains $p$~times~$+1$, $q$~times~$-1$,
  and
  $m$~subpairs~$(-1, +1)$, then there exist $\yy_1$, \pp,
  $\yy_p$ and $\zz_1$, \pp, $\zz_q$ in~$M$ such that
  $\xx_1^{e_1} \cdots \xx_n^{e_n} = \yy_1 \cdots \yy_p
  \zz_q\ii \cdots \zz_1\ii$ holds in~$G$, and the
  equality can be established using $m$~relations
  of~$\RX$. Then $\xx_1^{e_1} \cdots \xx_n^{e_n} = 1$ 
  in~$G$ implies $\yy_1 \cdots \yy_p = \zz_1
  \cdots \zz_q$ in~$G$, hence in~$M$. By
  Prop.~\ref{P:pres} again,
  the latter equality, if true, can be established using
  at most $(p+q)^2$ relations of~$\RX$. We have $m \le
  pq \le (p+q)^2/4$ and $p + q = n$, so, alltogether,
  we need $5n^2/4$ relations of~$\RX$ at most to
  establish $\xx_1^{e_1} \cdots \xx_n^{e_n} = 1$.
\end{proof}

So, for instance, all Garside groups satisfy a
quadratic isoperimetric inequality---as was already proved
in~\cite{Dfx}---and so do the groups~$G_i$ of
Example~\ref{X:maio}.

If we consider a Ore monoid~$M$ which is
\q-thin only, there is no clear result, as the
complexity of the group of units is involved. If the
group~$\UM$ is finite, then the presentation of
Prop.~\ref{P:pret} is finite, and it gives rise to a
quadratic isoperimetric inequality: however, in
this case, $M$ is thin, and Prop.~\ref{P:quis}
applies.

\section{Primitive elements}\label{S:atom}

The rest of the paper is centered on the
possible existence of an automatic structure for a
thin group of fractions, a strengthening of the
existence of a quadratic isoperimetric inequality. In
this section, we show that every thin monoid
admits a minimal spanning subset, whose elements will
be called primitive. Such primitive elements are
themselves connected with atoms, and we begin
with some observations about such elements
which extend earlier results of~\cite{Dfx, Dgk} so
as to allow the existence of nontrivial units.

\begin{defi}
  Assume that $M$ is a cancellative monoid. We say that an
  element~$\xx$ of~$M$ is an {\it atom} if $\xx$ is not
  a unit but $\xx = \yy \zz$ implies that either $\yy$
  or $\zz$ is a unit. The set of all atoms in~$M$ is
  denoted by~$\AM$.
\end{defi} 

\begin{lemm}\label{L:atcl}
  If $M$ is a cancellative monoid, then the
  atoms of~$M$ are closed under left and right
  multiplication by a unit, and, therefore,
  $\simeq$-saturated.
\end{lemm}

\begin{proof}
  Assume $\xx \in \AM$ and $\uu \in \UM$.
  First, $\xx \uu \in \UM$ would imply $\xx \in \UM$, so
  we have $\xx\uu \notin \UM$. Then assume $\xx\uu =
  \yy\zz$. We have $\xx = \yy(\zz
  \uu\ii)$, hence either $\yy$ or $\zz\uu\ii$ is a unit,
  and, in the latter case, so $\zz$. So $\xx\uu$ is an atom.
  The case of left multiplication by a unit is similar.
\end{proof}

\begin{defi}
  Assume that $M$ is a monoid. For~$\xx \in M$, we
  put $\n(\xx) = 0$ if $\xx$ is a unit, and
  \[
    \n(\xx) = \sup\{ n \, ; \, (\exists \xx_1, \pp, \xx_n
    \in \AM)(\xx = \xx_1 \cdots \xx_n)\}
  \]
  otherwise, if such a decomposition 
  exists and the involved supremum is finite.  
  Then we say that $M$ is {\it  \q-atomic} if
  $\n(\xx)$ exists for each~$\xx$ in~$M$. 
\end{defi}

Thus, $M$ is \q-atomic if and only if it is generated by
its atoms and its units, and, moreover, for
each~$\xx$ in~$M$, the maximal number of atoms
occurring in a decomposition of~$\xx$ is finite.
If there is no nontrivial unit in the monoid~$M$,
then $M$ is \q-atomic if and only if it is atomic in the
sense of~\cite{Dfx, Dgk}.

\begin{lemm}\label{L/norm}
  Assume that $M$ is a \q-atomic monoid. Then we have
  \begin{gather}
    \label{E:nor1}
    \n(\xx) = 0 \Leftrightarrow \xx \in \UM,\\
    \label{E:nor2}
    \n(\xx \yy) \ge \n(\xx) + \n(\yy),\\
    \label{E:nor3}
    \uu \in \UM \Rightarrow
    \n(\xx \uu) = \n(\uu \xx) = \n(\xx),
  \end{gather}
  for all~$\xx, \yy$ in~$M$. So $\xx \simeq
\xx'$ implies $\n(\xx) =
  \n(\xx')$.
\end{lemm}

\begin{proof}
  For~\eqref{E:nor1}, $\xx \in \UM$ implies $\n(\xx) = 0$
  by definition. Conversely, for $\xx \notin \UM$,
  the hypothesis that $\n(\xx)$ exists means that
  $\xx$ can be expressed as a finite product of atoms,
  so, by definition, $\n(\xx) \ge 1$ holds in this case.
  
  Assume now $\n(\xx) = p$ and $\n(\yy) = q$. For $p = q
  = 0$, both $\xx$ and $\yy$ are units, so is
  $\xx \yy$, and we have $\n(\xx\yy) = 0 = p + q$.
  For $p > 0$ and $q = 0$, $\yy$ is a unit, while
  $\xx$ admits a decomposition $\xx = \xx_1 \cdots
  \xx_p$ with $\xx_1$, \pp, $\xx_p \in \AM$. Then
  $\xx_p\yy$ is an atom, and we obtain $\xx\yy = \xx_1
  \cdots \xx_{p-1} (\xx_p\yy)$, hence $\n(\xx\yy) \ge p
  = p + q$. The argument is similar for $p = 0$ and $q >
  0$. Assume now $p > 0$ and $q > 0$. Then $\xx$
  and $\yy$ admit decompositions $\xx = \xx_1 \cdots
  \xx_p$, $\yy = \yy_1 \cdots \yy_q$ with $\xx_1$, 
  \pp, $\yy_q \in \AM$, and we deduce $\xx \yy = 
  \xx_1 \cdots \xx_p \yy_1 \cdots \yy_q$, hence $\n(\xx
  \yy) \ge p + q$. This shows \eqref{E:nor2}.
  
  Assume $\uu \in \UM$.
  Then \eqref{E:nor2} gives $\n(\xx \uu)
  \ge \n(\xx)$. Applying this with~$\xx\uu$ instead
  of~$\xx$ and $\uu\ii$ instead of~$\uu$, we obtain
  $\n({(\xx\uu)} \uu\ii) \ge \n(\xx\uu)$, \ie, 
  $\n(\xx) \ge \n(\xx \uu)$, whence $\n(\xx\uu) =
  \n(\xx)$. The argument is similar for~$\uu\xx$.
\end{proof}

\begin{prop}\label{P:atom}
  A monoid~$M$ is \q-atomic if and only if there
  exists a mapping $\l : M
  \rightarrow \NN$
  satisfying, for all~$\xx$, $\yy$ in~$M$,
  \begin{gather}
    \label{E:len1}
    \l(\xx) = 0 \Rightarrow \xx \in \UM,\\
    \label{E:len2}
    \l(\xx \yy) \ge \l(\xx) + \l(\yy).
  \end{gather}
\end{prop}

\begin{proof}
  Lemma~\ref{L/norm} shows that the mapping~$\n( )$
  satisfies \eqref{E:len1} and \eqref{E:len2} when $M$ is
  \q-atomic, so the condition is necessary.
  
  Conversely, assume that $\l$ is a mapping satisfying
  \eqref{E:len1} and \eqref{E:len2}. Assume $\xx \in
  M - \UM$, and let $\xx = \xx_1 \cdots \xx_p$ be a
  decomposition of~$\xx$ into non-invertible elements.
  By~\eqref{E:len1}, we have $\l(\xx_i) \ge 1$ for
  each~$i$, and, by~\eqref{E:len2}, 
  $\l(\xx) \ge \l(\xx_1) + \cdots + \l(\xx_p)$, hence
  $\l(\xx) \ge p$. So, the supremum~$n$ of the lengths
  of the decompositions of~$x$ into a product of 
  non-invertible elements satisfies $n \le l(\xx)$, and,
  therefore, it is finite. Now, let $\xx = \xx_1
  \cdots \xx_n$ be such a decomposition with maximal
  length. We claim that $\xx_1$, \pp, $\xx_n$ are
  atoms. Indeed, if $\xx_i$ is not an atom, it can be
  decomposed as $\xx_i = \xx'_i \xx''_i$ with
  neither~$\xx'_i$ nor~$\xx''_i$ in~$\UM$, and replacing
  $\xx_i$ with~$\xx'_i \xx''_i$ gives a decomposition
  of~$\xx$ of length~$n+1$. Hence every
  non-invertible element~$\xx$ of~$M$ is a product of
  at most $\l(\xx)$~atoms. This shows that $M$ is
  \q-atomic, and that $\n(\xx) \le \l(\xx)$ holds for
  every~$\xx$ in~$M$.
\end{proof}

\begin{exam}
  The monoids~$M_i$ of
  Example~\ref{X:main} all are \q-atomic, and
  even atomic as they contain no nontrivial unit: as the
  defining relations preserve the length, the
  latter induces a well defined mapping~$l$ of~$M_i$
  to~$\NN$ that satisfies~\eqref{E:len1} and~\eqref{E:len2}.
\end{exam}

The previous situation is general, as we have:

\begin{prop}\label{P:that}
  Every thin cancellative monoid is \q-atomic.
\end{prop}

\begin{proof}
  Assume that $M$ is a thin cancellative monoid, and
  $\SS$ is a finite spanning subset of~$M$. Let $\xx \in M$. As
  $\SS$ generates~$M$, $\xx$
  belongs to~$\SS^p$ for some~$p$. Now, let $\xx = \xx_1
  \cdots \xx_n$ be an arbitrary decomposition of~$\xx$
  with $\xx_1$, \pp, $\xx_n \in \SS - \UM$. By
  Prop.~\ref{P:clos}, $\SS^p$ spans~$M$, so the
  element~$\xx_i \cdots \xx_n$, which is a right divisor
  of~$\xx$, belongs to~$\SS^p$ as well. Assume
  $n \ge \card(\SS^p)$. Then there exist $i$, $j$ with
  $0 \le i < j \le n$ satisfying $\xx_i \cdots \xx_n = \xx_j
  \cdots \xx_n$, which implies $\xx_i \cdots \xx_{j-1} =
  1$ and  contradicts the assumption $\xx_i \notin \UM$.
  Thus we must have $n \le \card(\SS^p)
  \le \card(\SS)^p$. Let us define $\l(\xx)$ to be the
  maximal possible value of~$n$ in a decomposition as
  above. It is clear that $\l(\xx) = 0$ implies $\xx \in \UM$,
  and that $\l(\xx\yy) \ge \l(\xx) + \l(\yy)$ always holds, 
  as shows concatenating a maximal decomposition
  for~$\xx$ and a maximal decomposition for~$\yy$. So
  the mapping~$\l$ satisfies the
  conditions~\eqref{E:len1} and~\eqref{E:len2}, and, by
  Prop.~\ref{P:atom}, $M$ is \q-atomic.
\end{proof}

When we only assume that a finite \q-spanning set
exists, the situation is more complicated. However, we can
still recognize \q-atomicity as follows:

\begin{prop}
  Assume that $M$ is a \q-thin cancellative monoid.
  Then $M$ is \q-atomic if and only if the group~$\UM$
  is closed under conjugation
  in~$M$, in the sense that, if $\xx \uu = \uu' \xx$ holds,
  then $\uu \in \UM$ implies $\uu' \in \UM$.
\end{prop}
  
\begin{proof}
  Assume that $M$ is \q-atomic and we
  have $\xx \uu = \uu' \xx$ with $\uu \in \UM$.
  By Lemma~\ref{L/norm},
  we have $\n(\xx) = \n(\xx \uu) = \n(\uu' \xx) \ge
  \n(\uu') + \n(\xx)$, hence $\n(\uu') = 0$, and 
  $\uu' \in \UM$ by~\eqref{E:nor1}. So the condition
  is necessary.
  
  Conversely, assume that $\UM$ is closed under
  conjugation, and $\SS$ is a finite \q-spanning set
in~$M$.
  We adapt the argument of the proof of Prop.~\ref{P:that}.
  Let $\xx \in M$.  Then $\xx$ belongs to~$\SS^p \op
  \UM$ for some~$p$. Let $\xx = \xx_1 \cdots
  \xx_n \uu$ be any decomposition of~$\xx$ with
  with $\xx_1$, \pp, $\xx_n \in \SS - \UM$
  and $\uu \in \UM$. The set $\SS^p \op \UM$ 
  spans~$M$, so, by Prop.~\ref{P:clos}, $\xx_i \cdots
  \xx_n$, which is a right divisor of~$\xx$, belongs
  to~$\SS^p \op\UM$ for each~$i$. If $n \ge
  \card(\SS^p)$ holds, there exist $i$, $j$ with
  $0 \le i < j \le n$ satisfying $\xx_i \cdots \xx_n \simeq
  \xx_j \cdots \xx_n$, so we have
  $(\xx_i \cdots \xx_{j-1})(\xx_j \cdots \xx_n) = 
  (\xx_j \cdots \xx_n) \uu$ for some unit~$\uu$. This, by
  hypothesis, implies $\xx_i \cdots \xx_{j-1} \in \UM$,
  contradicting the hypothesis $\xx_i \notin \UM$. So
  we must have $n \le \card(\SS^p) \le \card(\SS)^p$. If
  we define $\l(\xx)$ to be the
  maximal possible value of~$n$ in a decomposition as
  above, then $\l$
  satisfies~\eqref{E:len1} and~\eqref{E:len2},
  and $M$ is \q-atomic.
\end{proof}

If $M$ is a \q-atomic cancellative monoid, then the
relation~$\div$ is a strict partial ordering
with no infinite descending chain (and so is its right
counterpart). Indeed, by~\eqref{E:nor2}, $\xx \div
\yy$ implies $\n(\xx) < \n(\yy)$, so $\div$ may
admit no cycle, hence it is a strict partial ordering,
and it admits no infinite descending chain since
$(\NN, <)$ does. In such a framework, we can
introduce the notion of a minimal common (right)
multiple (``mcm''), which extends the notion of a
least common multiple (``lcm''):

\begin{defi}
  Assume that $M$ is a monoid. For $\xx$, $\yy \in M$, we
  say that $\zz$ is an {\it mcm} of~$\xx$ and~$\yy$
  if both $\xx\dive\zz$ and $\yy\dive\zz$ hold, but
  $\xx\dive\zz'$ and $\yy\dive\zz'$ hold for no
  proper divisor~$\zz'$ of~$\zz$. We denote
  by~$\CC(\xx, \yy)$ the set of all elements~$\yy'$
  such that $\xx\yy'$ is an mcm of~$\xx$ and~$\yy$
  (if any).
\end{defi}

An mcm is like a lcm, except that we require no
uniqueness. For instance, in the monoid~$M_1$ of
Example~\ref{X:main}, the elements~$\a$ and~$\b$
admit two mcm's, namely $\a^2$ and $\a\b$, but they
admit no lcm, as we have neither $\a^2 \dive
\a\b$ nor $\a\b \dive \a^2$.

\begin{lemm}\label{L:rmcm}
  Assume that $M$ is a \q-atomic cancellative monoid.
  Then, for all~$\xx$, $\yy$ in~$M$, every common
  multiple of~$\xx$ and~$\yy$ (if any) is a multiple
  of some mcm of~$\xx$ and~$\yy$.
\end{lemm}

\begin{proof}
  Assume $\xx \dive \zz$ and $\yy \dive \zz$. Let $Z
  = \{\zz' \dive \zz \, ; \, \xx \dive \zz' , \yy \dive \zz'
  \}$. Then any element~$\zz'$ of~$Z$ such
  that $\n(\zz')$ has the least possible value is an
  mcm of~$\xx$ and~$\yy$.
\end{proof}

We are now ready to introduce the notion of a primitive
element:

\begin{defi}
  Assume that $M$ is a \q-atomic cancellative monoid.
  We say that an element~$\xx$ of~$M$ is {\it primitive}
  if $\xx$ belongs to the smallest subset~$\SS$ of~$M$
  that contains the atoms and is such that, if $\xx$ and
  $\yy$ belong to~$\SS$, so does $\yy'$ whenever
  $\xx\yy'$ is an mcm of~$\xx$ and~$\yy$. The set of
  all primitive elements of~$M$ is denoted by~$\PM$.
\end{defi}

  In other words, $\PM$ is the closure of~$\AM$
  under operation~$\CC$.

\begin{exam}
  If any two elements admitting a common multiple
  admit  a lcm, the set $\CC(\xx, \yy)$ is either
  empty, or it consists of a single element $\xx \dr \yy$,
  so the primitive elements are the closure of the atoms
  under operation~$\dr$.
  
  Consider now the monoid~$M_1$ of
  Example~\ref{X:main}. There are two atoms, namely $\a$
  and~$\b$. We observed above that $\a$ and $\b$
  admit two right mcm's, namely $\a^2$ and $\a\b$,
  so $\CC(\a, \b)$ consist of the two elements~$\a$
  and~$\b$, and so does $\CC(\b, \a)$. It follows that
  the closure of~$\AM$ under~$\CC$ is the set $\{1,
  \a, \b\}$, \ie, there are three primitive elements
  in~$M_1$.
  
  The reader can easily check that there are four primitive
  elements in the monoids~$M_2$ and~$M_3$, namely
  $1$ and the atoms $\a$, $\b$, and~$\c$.
\end{exam}

\begin{lemm}\label{L:prsa}
  Assume that $M$ is a \q-atomic cancellative monoid.
  Then the set~$\PM$ is closed under right
  multiplication by a unit, and, therefore,
  it is $\simeq$-saturated.
\end{lemm}

\begin{proof}
  Assume $\xx \in \PM$, and $\uu \in \UM$. If
  $\xx$ is an atom, then $\xx \uu$ is an atom as well, so
  it is primitive. Otherwise, there exist~$\yy$, $\zz$
  in~$\PM$ such that $\yy\xx$ is an mcm
  of~$\yy$ and~$\zz$. In this case, $\yy\xx\uu$ is also an mcm of~$\yy$ and~$\zz$, so $\xx\uu$ also
  belongs to~$\CC(\yy, \zz)$, and, therefore, it is
  primitive. So, we have $\PM \op \UM \ince \PM$, and,
  therefore, $\PM$ is $\simeq$-saturated.
\end{proof}

We can now prove:

\begin{prop}\label{P:micl}
  Assume that $M$ is a \q-atomic cancellative monoid.
  Then $\PM$ spans~$M$, and every spanning subset
  of~$M$ includes~$\PM$.
\end{prop}

\begin{proof}
  First, $\PM$ is $\simeq$-saturated by
  Lemma~\ref{L:prsa}. Then, if $\xx$ is an atom of~$M$
  and $\uu$ is a unit, $\xx \uu$ is an mcm of~$\xx$
  and~$\xx$, and, therefore, $\uu$ is primitive. Thus
  $\PM$ includes $\AM$ and $\UM$, and, therefore, it
  generates~$M$. Next, assume $\xx \dive \zz$ and $\yy
  \dive \zz$ with $\xx$, $\yy \in \PM$. By
  Lemma~\ref{L:rmcm}, there exist
  $\xx'$, $\yy'$ such that $\xx\yy' = \yy \xx' \dive \zz$
  holds and $\xx\yy'$ is an mcm of~$\xx$ and~$\yy$.
  This implies $\xx' \in \PM$ and $\yy' \in \PM$ by
  definition. Hence $\PM$ satisfies Condition~\eqref{E:coni}, \ie, it spans~$M$.
  
  Let $\SS$ be an arbitrary spanning subset of~$M$.
  As $\SS$ generates~$M$, it necessarily includes
  $\AM$. Then, $\SS$ has to be closed under~$\CC$. 
  Indeed, assume that $\xx$, $\yy$ lie in~$\SS$, and
  $\xx\yy''$ is an mcm of~$\xx$ and~$\yy$, say
  $\xx\yy'' = \yy\xx''$. As $\SS$ spans~$M$,
  there exist $\yy'$ in~$\SS$
  and $\zz$ satisfying $\yy'' = \yy' \zz$.
  The hypothesis that $\xx\yy''$ is an mcm implies
  $\n(\xx\yy') = \n(\xx\yy'')$, hence $\n(\zz) = 0$,
  and $\yy' \simeq \yy''$. As $\SS$ is
  $\simeq$-saturated by definition, we deduce $\yy''
  \in \SS$. So, $\SS$
  includes~$\AM$ and it is closed under~$\CC$, hence it
  includes the closure~$\PM$ of~$\AM$ under~$\CC$.
\end{proof}

\begin{coro}
  Assume that $M$ is a \q-atomic cancellative monoid.
  Then $M$ is thin (\resp \q-thin) if and only the
  set~$\PM$ is finite (\resp \q-finite).
\end{coro}

Once we know that primitive elements span~$M$,
we can apply Prop.~\ref{P:clos} and we deduce that, if
$M$ is a \q-atomic cancellative monoid, then the
set~$\PM$ is closed under left multiplication by a
unit, and every right divisor of a primitive element is
primitive.

Another consequence of Prop.~\ref{P:micl} is
that, if $M$ is a \q-atomic cancellative monoid,
and $\SS$ is a $\simeq$-selector through~$\PM$,
then $\SS$ is a minimal \q-spanning subset of~$M$, 
and $\SS \cap \AM$ is a $\simeq$-selector 
through~$\AM$. Indeed, by construction, $\SSU$ is
equal to~$\PM$, so $\SS$ \q-spans~$M$.
If $\SS'$ is a proper subset of~$\SS$, $\SS' \op \UM$
is a proper subset of~$\PM$, so it cannot span~$M$,
and $\SS'$ cannot \q-span~$M$. Finally, every
atom~$\xx$ is primitive, so it belongs to~$\SSU$,
and, therefore, $\xx$ is
$\simeq$-equivalent to one element of~$\SS$. 

\section{Simple elements}\label{S:simp}

A crucial feature in Garside's and Thurston's analysis of
the braid monoids and its subsequent extensions is
the existence of a finite subset that is closed both
under lcm and right divisor: in the current
framework, this means that there exists a finite
spanning subset~$\SS$ that is closed under lcm, \ie,
the lcm of two elements of~$\SS$ belongs to~$\SS$.
The least such set~$\SS$ happens to be the closure of
primitive elements under lcm, and its elements, called
minimal in~\cite{Cha, Chb}, or simple in~\cite{Dgk},
play a prominent r\^ole. In particular, there exists a
maximal simple element~$\D$ which enjoys most of
the properties of Garside's fundamental
braids~$\D_n$ \cite{Gar}.

So, in the current approach, a natural idea would be to 
look for finite spanning subsets closed under mcm. 
Unfortunately, when least
common multiples do not exist, more precisely when
common  multiples exist but some elements admit at least
two non $\simeq$-equivalent mcm's, no such set may
exist:

\begin{prop}\label{P:node}
  Assume that $M$ is a \q-atomic cancellative monoid,
  any two elements of~$M$ admit a common
  multiple, and $\SS$ is a finite spanning subset of~$M$ that
  is closed under mcm. Assume in addition
  that $\xx \dive \yy \in \SS$ implies $\xx \in \SS$. Then 
  any two elements of~$M$ admit a lcm.
\end{prop}

\begin{proof}
  As $\SS$ is finite and closed under right mcm,
  there exists~$\D$ in~$\SS$ such that $\xx \dive \D$
  holds for every~$\xx$ in~$\SS$, \ie, there
  exists~$\xx\dd$ satisfying $\xx \xx\dd = \D$; as $M$ is
  left cancellative, $\xx\dd$ is unique, and, as $\SS$ 
  spans~$M$, every right divisor of an
  element of~$\SS$ belongs to~$\SS$, so $\xx\dd$ belongs
  to~$\SS$. The mapping $\xx \mapsto \xx\dd$ is injective,
  and, therefore, it is a permutation of~$\SS$. Assume
  $\xx$, $\yy \in \SS$, and let $\xx\yy'$ and $\xx\yy''$
be two
  right mcm's of~$\xx$ and $\yy$. By hypothesis, $\xx\yy'$
  belongs to~$\SS$, so, by the previous remark, there
  exists~$\zz$ in~$\SS$ satisfying $\xx\yy' = \zz\dd$, \ie,
  $\zz \xx \yy' = \D$. By hypothesis, $\SS$ is closed
  under left divisors, so $\zz\xx$ and,
  similarly, $\zz\yy$ belong
  to~$\SS$, and so does their right mcm $\zz\xx\yy''$. So
  we must have $\zz\xx\yy'' \dive \D$, hence $\yy'' \dive
  \yy'$. In other words, $\xx\yy'$ is a lcm of~$\xx$
  and~$\yy$. Finally, as $\SS$ generates~$M$, the
  existence of a lcm for each pair of elements of~$\SS$
  inductively implies the existence of a lcm for each
  pair of elements of~$M$.
\end{proof}

Thus, we must find a more subtle definition. The
following one is convenient, in the sense that
it will prove appropriate for the construction of a normal
form.

\begin{defi}
  Assume that $M$ is a cancellative monoid, and 
  $\SS$ \q-spans~$M$. We say that an
  element~$\xx$ of~$M$ is {\it $\SS$-simple} if $\yy
  \div \xx$ implies $\Dl(\yy) \cap \SS \not= 
  \Dl(\xx) \cap \SS$. If $M$ is \q-atomic, we say {\it
  simple} for $\PM$-simple.
\end{defi}

The elements of~$\SS$ always are $\SS$-simple. Indeed,
for $\xx \in \SS$, we have $\xx \in
\Dl(\xx) \cap \SS$, but $\yy \div \xx$ implies $\xx
\not\div \yy$, \ie, $\xx \notin \Dl(\yy) \cap \SS$.
So, in particular, a primitive element is always simple.
By definition, an element is $\SS$-simple if and only if it is
an mcm of its divisors lying
in~$\SS$. In particular, in the Gaussian case, an element
is $\SS$-simple if and only if it is the lcm of its
divisors lying in~$\SS$, and, therefore, a $\SS$-simple
element~$\xx$ is completely determined by the
set~$\Dl(\xx) \cap \SS$. This need not be true in the
general case.
  
\begin{exam}\label{X:simp}
  Let $M$ be a free commutative monoid based
  on~$\{\a_1, \pp, \a_n\}$. Then the atoms
  of~$M$ are $\a_1$, \pp, $\a_n$, there are $n+1$
  primitive elements, namely $1$ and the atoms, and there
  are $2^n$~simple elements, namely the elements
  $\a_I = \prod_{i \in I}\a_i$ for $I \ince \{1, \pp, n\}$.
  Indeed, $\a_i \dive \a_I$ is equivalent to~$i \in I$
  in this case, and, for every~$\xx$, the element~$\a_I$
  with $I = \{i \,;\, \a_i \dive x\}$ is a divisor
  of~$\xx$
  with the same divisors in~$\PM$, so no element not of
  the form~$\a_I$ may be simple.

  On the other hand, there are three primitive elements,
  namely $1$, $\a$, and~$\b$ in the monoid~$M_1$ of
  Example~\ref{X:main}. These elements are simple, and
  there are two more simple elements, namely $\a^2$ and
  $\a\b$. Here, we have $\Dl(\a^2) \cap \PM
  = \Dl(\a\b) \cap \PM = \PM$, which gives an
  example where a simple element is not
  determined by the family of its primitive divisors.
  
  Similarly, there are seven simple elements in~$M_3$,
  namely the four primitive elements $1$, $\a$, $\b$,
  $\c$, and, in addition, the three elements $\a\b$,
  $\b\a$, and $\b^2$: the sets of primitive
divisors of the latter
  elements are $\{\a, \b\}$, $\{\b, \c\}$, and $\{\a, \b, 
  \c\}$ respectively, so, here, a simple element
  happens to be determined by its primitive
  divisors (although the monoid admits no lcm).
\end{exam}

In the Gaussian case, \ie, when least common
multiples exist, the current definition of a simple
element is equivalent to that of~\cite{Dgk}. In
particular, in the thin case, the simple elements
of~$M$ have a natural characterization extending
that given for a free commutative monoid in
Example~\ref{X:simp}. 
  
\begin{prop}\label{P:gasi}
  Assume that $M$ is a thin Gaussian monoid, \ie, a 
  Garside monoid. Then the
  simple elements of~$M$ are finite in number, and they
  coincide with the divisors of the lcm~$\D$ of~$\PM$.
\end{prop}

\begin{proof}
  Let $\{\xx_i \, ; \, i = 1, \pp, n \}$ be an enumeration
  of~$\PM$. For $I \ince \{1, \pp, n\}$, let $\xx_I$ be the
  lcm of the~$\xx_i$'s with $i \in I$. Then $\xx_I$ is
  simple, and, conversely, every simple element must be of
  this form. Let $\D$ be the lcm of~$\PM$. Then, 
  by construction, every simple element~$\xx_I$ is 
  a divisor of~$\D$. The computation rules for lcm's
  then imply that simple elements span~$M$
  \cite{Dgk}, and, as a consequence, that every
  divisor of~$\D$ is simple. 
\end{proof}

It is well known that, if $M$ is a spherical Artin monoid,
then the simple elements are in one-to-one correspondence
with the elements of the associated finite Coxeter group
\cite{BrS, Dlg}: for instance, the $n!$~simple elements in
the braid monoid~$B_n^+$ are in one-to-one
correspondence with the permutations of $n$~objects.
More generally, it is shown in~\cite{Dgk} that the simple
elements of a Garside monoid make a finite lattice
with a unique maximal element, the lcm~$\D$ of
the primitive element. As shows the case of the
monoid~$M_1$, such a property need not be true in the
general case.

For future use, we gather now some general
results about $\SS$-simple elements.

\begin{lemm}\label{L:sisa}
  Assume that $M$ is a cancellative monoid,
  and $\SS$ \q-spans~$M$.
  
  \noindent (i) An element of~$M$ is $\SS$-simple if
  and only if it is $\SSU$-simple.
  
  \noindent (ii) The set of all $\SS$-simple elements
  is closed
  under left and right multiplication by a
  unit, and, therefore, it is $\simeq$-saturated.
\end{lemm}

\begin{proof}
  (i) Assume that $\xx$ is $\SS$-simple, and $\yy \div
  \xx$ holds. By definition, we have $\Dl(\yy) \cap \SS
  \not= \Dl(\xx) \cap \SS$, so, a fortiori, $\Dl(\yy) \cap 
  \SSU \not= \Dl(\xx) \cap \SSU$, and $\xx$ is
  $\SSU$-simple. Conversely, assume that $\xx$ is
  $\SSU$-simple, and $\yy \div \xx$ holds. Then we
  have $\zz \dive \xx$ and $\zz \not\dive \yy$ for some
  $\zz$ in~$\SSU$. By definition, we have $\zz = \zz'
  \uu$ for some~$\zz'$ in~$\SS$ and $\uu$
in~$\UM$.
  Then $\zz' \dive \xx$ and $\zz' \not\dive \yy$ hold,
  and $\xx$ is $\SS$-simple.
  
  (ii) Assume that $\xx$ is $\SS$-simple and $\uu$ is
  a unit. Then we have $\Dl(\xx \uu) = \Dl(\xx)$, and
  $\yy \div \xx \uu$ is equivalent to~$\yy \div \xx$.
  Hence $\yy \div \xx$ implies $\Dl(\yy) \cap \SS \not=
  \Dl(\xx \uu)\cap \SS$, and $\xx\uu$ is $\SS$-simple.
  The argument is similar for left multiplication by~$\uu$,
  as $\Dl(\uu \xx) = \Dl(\xx)$ holds.
\end{proof}

In the general case, as shows the example
of~$M_1$, simple elements need not span the
monoid. However, we still have the following closure 
property:

\begin{lemm}\label{L:ridi}
  Assume that $M$ is a cancellative monoid,
  and $\SS$ \q-spans~$M$.
  Then every right divisor of a $\SS$-simple element
  of~$M$ is $\SS$-simple.
\end{lemm}

\begin{proof}
  Assume that $\xx\yy$ is $\SS$-simple. We wish to show
  that $\yy$ is $\SS$-simple. Assume $\yy' \div \yy$.
  Then we have $\xx\yy' \div \xx\yy$, so, by definition,
  there exists $\zz$ in~$\SS$ satisfying $\zz \dive \xx
  \yy$ and $\zz \not\dive \xx\yy'$. As $\SSU$ spans~$M$,
  there exists $\zz'$ in~$\SSU$ and $\xx'$
  in~$M$ satisfying $\xx\zz' = \zz\xx' \dive \xx\yy$,
  hence $\zz' \dive \yy$. Thus, $\zz'$ belongs
  to~$\Dl(\yy) \cap \SSU$. On
  the other hand, $\zz' \dive \yy'$ would
  imply $\zz\xx' = \xx\zz' \dive \xx\yy'$, hence $\zz
  \dive \xx\yy'$, contradicting the hypothesis. So $\zz'
  \dive \yy'$ is impossible. We have
  $\Dl(\yy') \cap \SSU \not= \Dl(\yy) \cap \SSU$, so
  $\yy$ is $\SSU$-simple, hence $\SS$-simple.
\end{proof}

Proposition~\ref{P:gasi} implies that, in the Gaussian case,
there exist at most $2^n$~simple elements when there are
$n$~primitive elements, a bound which we have seen is
nearly reached in the case of a free commutative monoid.
The result extends to the general case as follows:

\begin{prop}\label{P:fisi}
  Assume that $M$ is a cancellative
  monoid, and $\SS$ is a \q-spanning subset of~$M$
  with $n$~elements. Then every $\SS$-simple
  element belongs to $\SS^n \op \UM$, and, therefore,
  there are at most $n^n$ $\simeq$-equivalence classes of
  $\SS$-simple elements in~$M$. 
\end{prop}

\begin{proof}
  Assume that $\xx$ is $\SS$-simple, and let $\xx_1$,
  \pp, $\xx_p$ be an enumeration of~$\Dl(\xx) \cap\SS$.
  By Proposition~\ref{P:powf}(iii), there exists~$\xx'$
  in~$\SS^p$
  such that $\xx_i \dive \xx'$ holds for $1 \le i \le p$.
  Then,  we have $\Dl(\xx') \cap \SS
  \supseteq \Dl(\xx) \cap \SS$, hence $\Dl(\xx') \cap \SS
  = \Dl(\xx) \cap \SS$.
  By definition of a $\SS$-simple element, $\xx' \div
  \xx$ is impossible, so $\xx' \simeq \xx$ is the only
  possibility, which shows that $\xx$ belongs to~$\SS^p
  \op \UM$, and, therefore, to~$\SS^n \op \UM$ since $1$
  is primitive and $p \le n$ holds.
\end{proof}

\begin{coro}\label{C:fisi}
  If $M$ is a thin cancellative monoid, then the
  simple elements of~$M$ are finite in number. More
  precisely, if $M$ contains $n$~primitive elements, it
  contains at most $n^n$~simple elements.
\end{coro}

If $M$ is a \q-atomic cancellative monoid and $S$ is a
$\simeq$-selector through simple elements, then $S
\cap \PM$ is a $\simeq$-selector through~$\PM$,
and $S \cap \AM$ is a $\simeq$-selector
through~$\AM$. Conversely, every $\simeq$-selector
through~$\AM$ can be extended into a selector
through~$\PM$, and, then, through simple elements.

Finally, let us observe that simple elements, as atoms
and primitive elements, are defined intrinsically,
and, therefore, they are preserved by automorphisms:

\begin{prop}\label{P:prau}
  Assume that $M$ is a \q-atomic cancellative monoid,
  and $\phi$
  is an automorphism of~$M$. Then $\phi$ globally
  preserves $\AM$, $\PM$, and the set of all simple
  elements in~$M$.
\end{prop}

\begin{proof}
  As $\phi$ maps units to units, then it maps
  non-atoms to non-atoms, and, therefore, it
  maps atoms to atoms. Then, it maps every mcm
  of two elements to an mcm of their images, and,
  therefore, it maps every primitive element to a 
  primitive element. Finally, $\phi$ preserves the
  relations~$\dive$ and~$\div$, so it maps simple 
  elements to simple elements.
\end{proof}

\section{Normal forms}\label{S:nofo}

The main interest of simple elements in the Gaussian
case, \ie, when least common multiples exist, 
is that they can be used to construct good
normal forms. In particular, the greedy normal
form originally defined for the braid monoids
\cite{Adj, Thu, ElM, Eps} extends to every Gaussian
monoid, and, subsequently, to the corresponding
group of fractions \cite{Dgk}. The principle is that, for
~$\xx \not= 1$ in the considered monoid~$M$, there
exists a maximal simple divisor~$\xx_1$ of~$\xx$,
namely the gcd of~$\xx$ and the maximal simple
element~$\D$, so we can write $\xx = \xx_1 \xx'$,
and, applying the process to~$\xx'$, we inductively
obtain a decomposition $\xx = \xx_1 \xx_2 \cdots$ in
terms of simple elements. This decomposition enjoys
good properties, and, in particular, it gives rise to a
bi-automatic structure on the associated group of
fractions.

A crucial technical point in the above construction is that
simple elements happen to span the monoid, in the
Gaussian case. We shall see now that a similar
construction is still possible in the general case when
we start with an arbitrary spanning set~$\SS$ and use
the derived notion of a $\SS$-simple element. The
price to pay for the generalization is that a given
element possibly may have more than one normal
decomposition, but, this fact excepted, the results
remain similar, and the proofs are extremely easy.

As in the Gaussian case, we start from the fact that, for
every element~$\xx$, there exists a maximal $\SS$-simple
divisor of~$\xx$:

\begin{lemm}\label{L:nor0}
  Assume that $M$ is a \q-atomic cancellative monoid,
  $\SS$ \q-spans~$M$, and $\xx_0$
  is a $\SS$-simple element. Then, for every~$\xx$
  in~$M$, there exists a $\SS$-simple divisor~$\xx_1$ of~$\xx$ satisfying
  $\Dl(\xx) \cap \SS = \Dl(\xx_1) \cap \SS$; moreover, 
  we may assume $\xx_0 \dive \xx_1$ whenever $\xx_0
  \dive \xx$ holds.

  If $M$ is Gaussian, then $\xx_1$ is
  a lcm of $\Dl(\xx) \cap \SS$, and, so, it is 
  unique.
\end{lemm}

\begin{proof}
  Assume $\xx_0 \dive \xx$. Let $\YY$ be the set of all
  $\SS$-simple elements~$\yy$ satisfying $\xx_0 \dive \yy
  \dive \xx$, and let $\xx_1$ be an element of~$\YY$
  such that $\n(\xx_1)$ has the maximal possible value:
  such an element exists since $\yy \in \YY$ implies
  $\n(\yy) \le \n(\xx)$. Write $\xx = \xx_1 \xx''$.
  Assume $\zz \in \Dl(\xx) \cap \SS$. As $\SS$ \q-spans~$M$, there must exist $\zz'$
  in~$\SS$, and $\xx'_1$ in~$M$ satisfying
  $\xx_1 \zz' = \xx'_1 \zz' \dive \xx$. So we have
  $\xx_0 \dive \xx_1 \dive \xx_1 \zz' \dive \xx$.
  Moreover, provided
  $\zz'$ has been chosen so that $\n(\xx_1 \zz')$ has the
  least possible value, no proper divisor of~$\xx_1 \zz'$
  is a multiple of~$\zz$, which implies that $\xx_1\zz'$
  is $\SS$-simple, and, therefore, it belongs to~$\YY$.
  The definition of~$\xx_1$ then implies $\n(\xx_1\zz') = 
  \n(\xx_1)$, hence $\zz' \in \UM$, and then $\zz \dive
  \xx_1$. So we have $\Dl(\xx_1) \cap \SS =
  \Dl(\xx) \cap \SS$. Take $\xx_0 = 1$ for the
  general result.
  
  In the Gaussian case, the lcm of~$\Dl(\xx) \cap
  \SS$ is a $\SS$-simple element satisfying the
  requirements, and it divides every other element
  satisfying them, so it must be the only solution.
\end{proof}

\begin{defi}\label{D:nof1}
  Assume that $M$ is a cancellative monoid, 
  and $\SS$ \q-spans~$M$.
  We say that a sequence $(\xx_1, \pp, \xx_n)$ in~$M$ is
  {\it $\SS$-prenormal} if, for each~$i$, we have 
  $\Dl(\xx_i) \cap \SS = \Dl(\xx_i \cdots \xx_n) \cap
  \SS$. We say that $(\xx_1, \pp, \xx_n)$ is
  {\it $\SS$-normal} if it is $\SS$-prenormal, and, in
  addition, each factor~$\xx_i$ is $\SS$-simple.
  If $M$ is \q-atomic, we say (pre)normal for
  $\PM$-(pre)normal.
\end{defi}

Say that a sequence $(\xx_1, \pp, \xx_n)$ is a
decomposition for~$\xx$ if
$\xx = \xx_1 \cdots \xx_n$ holds. Iterating
Lemma~\ref{L:nor0}, we find:

\begin{prop}\label{P:norl}
  Assume that $M$ is a \q-atomic cancellative monoid, 
  $\SS$ \q-spans~$M$, and $\xx_0$
  is a $\SS$-simple element. Then every element~$\xx$
  of~$M$ satisfying $\xx_0 \dive \xx$ admits a $\SS$-normal decomposition $(\xx_1, \pp, \xx_n)$
  with $\xx_0 \dive \xx_1$.
\end{prop}

\begin{proof}
 Lemma~\ref{L:nor0} gives a $\SS$-simple
 element~$\xx_1$ satisfying both $\xx_0 \dive
\xx_1$ and
 $\Dl(\xx_1) \cap \SS = \Dl(\xx) \cap \SS$.
 Write $\xx = \xx_1 \xx'$. If $\xx'$ is a unit, then $\xx$
 is $\SS$-simple by Lemma~\ref{L:sisa}(ii), and we are
 done. Otherwise, we have $\n(\xx') < \n(\xx)$:
  inductively, we find a $\SS$-normal decomposition
  $(\xx_2, \pp, \xx_n)$ of~$\xx'$, and concatenating
  $\xx_1$ with the latter gives a $\SS$-normal
  decomposition of~$\xx$.
\end{proof}

By Lemma~\ref{L:nor0}, the normal form of
Prop.~\ref{P:norl} is unique in the Gaussian case, and it
coincides with the greedy normal form
of~\cite{Eps, Dgk}. More generally, the $\SS$-normal
form is unique whenever distinct $\SS$-normal
elements never admit the same divisors in~$\SS$.
Now, the latter condition need not be true, and the
normal form need not be unique in general.

\begin{exam}\label{X:norm}
  Consider once again the monoid~$M_1$ of
  Example~\ref{X:main}. Then $\a^2$ and $\a\b$ are
  simple elements, and $(\a^2, \a^2)$ and $(\a\b,
  \a\b)$ are two normal decompositions for~$\a^4$
  in~$M$. On the other hand, we observed that, in the
  case of~$M_3$, the simple elements are uniquely
  determined by their primitive divisors. So, in this case,
  the normal decomposition is unique.
\end{exam}

When nontrivial units exist, we can replace the
family of all $\SS$-simple elements by a
$\simeq$-selector, at the expense of keeping a unit at the
end of the decomposition. 

\begin{coro}
  Assume that $M$ is a \q-atomic cancellative monoid, 
  $\SS$ \q-spans~$M$, and $\S$ is a
  $\simeq$-selector through $\SS$-simple elements
  in~$M$. Then every element~$\xx$
  of~$M$ admits a decomposition 
  $(\xx_1, \pp, \xx_n, \uu)$ with $\xx_1$, 
  \pp, $\xx_n \in \S$, $\uu \in \UM$ and $\xx_i \cov(\SS)
  \xx_{i+1}$ for each~$i$.
\end{coro}

\begin{proof}
  Let $\xx \in M$. By
  Prop.~\ref{P:norl}, $\xx$ admits a $\SS$-normal
  decomposition $(\xx'_1, \pp,  \xx'_n)$, where $\xx'_1$, 
  \pp, $\xx'_n$ are $\SS$-simple and $\xx'_i \cov(\SS)
  \xx'_{i+1}$ holds for every~$i$. Using
  \eqref{E:poqu}, we find $\xx_1$, \pp, $\xx_n$ in~$\S$
  and $\uu$ in~$\UM$
  satisfying $\xx_i \simeq \xx'_i$ for every~$i$ and
  $\xx_1 \cdots \xx_n \uu = \xx'_1 \cdots \xx'_n = \xx$.
  By Lemma~\ref{L:exco},
  $\xx'_i \cov(\SS) \xx'_{i+1}$ implies 
  $\xx_i \cov(\SS) \xx_{i+1}$.
\end{proof}

The interest of the current construction lies in that
$\SS$-normal sequences admit a purely {\it local}
characterization. 

\begin{defi}
  Assume that $M$ is a monoid and $\SS$ is a subset 
  of~$M$. For $\xx$, $\yy \in M$, we say
  that $\xx$ {\it covers}~$\yy$ w.r.t.{ }$\SS$, denoted
  $\xx  \cov(\SS) \yy$, if $\Dl(\xx\yy) \cap \SS
  = \Dl(\xx) \cap \SS$ holds, \ie, if every element
  of~$\SS$ dividing~$\xx\yy$ already divides~$\xx$.
  If $M$ is \q-atomic, we write $\cov()$ for
  $\cov(\PM)$.
\end{defi}

\begin{lemm}\label{L:exco}
  Assume that $M$ is a cancellative monoid and $\SS$ is
  a subset of~$M$. Then the relations~$\cov(\SS)$ and
  $\cov(\SSU)$ coincide.
\end{lemm}

\begin{proof}
  It is obvious that $\xx \cov(\SSU) \yy$ implies
  $\xx \cov(\SS) \yy$. Conversely, assume $\xx
  \cov(\SS) \yy$ and $\zz' \simeq \zz \in \SS$. Then
  $\zz' \dive \xx\yy$ (\resp $\zz' \dive
  \xx$) is equivalent to $\zz \dive \xx\yy$ (\resp
  $\zz \dive \xx$). So $\zz' \dive \xx\yy$ implies
  $\zz \dive \xx\yy$, hence $\zz\dive \xx$, hence
  $\zz' \dive \xx$, and we have $\xx \cov(\SSU) \yy$. 
\end{proof}

All properties of the normal form relie on the following
basic observations:

\begin{lemm}\label{L:cov1}
  Assume that $M$ is a cancellative monoid, and 
  $\SS$ \q-spans~$M$. Then, for $\xx$, $\yy$,
  $\zz $ in~$M$:
  
  \noindent (i) The relation $\yy \cov(\SS) \zz$ implies
  $\xx \yy  \cov(\SS) \zz$; 
  
  \noindent (ii) The conjunction of $\xx \cov(\SS) \yy$
  and $\yy \cov(\SS) \zz$ implies $\xx \cov(\SS) \yy\zz$.
\end{lemm}

\begin{figure}
  \includegraphics{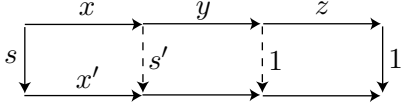}
   \caption{Covering relation}\label{F:cove}
\end{figure}

\begin{proof}
  (i) (Fig.~\ref{F:cove}) Assume $\ss \in \SS$ and $\ss
  \dive \xx \yy 
  \zz$. As $\SSU$ spans~$M$, there
  exist $\ss'$ in~$\SSU$, and $\xx'$ in~$M$
  satisfying $\ss \xx' = \xx \ss' \dive \xx \yy
  \zz$, hence $\ss' \dive \yy \zz$. By
  Lemma~\ref{L:exco}, $\yy \cov(\SS) 
  \zz$ implies $\yy \cov(\SSU) \zz$, so $\ss' \dive \yy \zz$
  implies $\ss' \dive \yy$, and, therefore, we have
  $\ss \dive \xx \ss' \dive \xx \yy$.
 
  (ii) Assume $\ss \in \SS$ and $\ss \dive \xx \yy \zz$. 
  By~(i), $\yy \cov(\SS) \zz$ implies  $\xx \yy \cov(\SS)
  \zz$, so we deduce $\ss \dive \xx \yy$. Then the
  hypothesis $\xx \cov(\SS) \yy$ implies $\ss \dive \xx$.
\end{proof}

We can know establish the expected local characterization
of normal sequences, a necessary first step toward a
possible automatic structure:

\begin{prop}\label{P:loca}
  Assume that $M$ is a cancellative monoid, and $\SS$ 
  \q-spans $M$. Then a sequence
  $(\xx_1, \pp,  \xx_n)$ in~$M^n$ is $\SS$-prenormal
  if and only if $\xx_i \cov(\SS) \xx_{i+1}$ holds for
  each~$i$.
\end{prop}

\begin{proof}
  By definition, the sequence~$(\xx_1, \pp, \xx_n)$ is
  $\SS$-prenormal if and only if $\xx_i \cov(\SS)
  \xx_{i+1} \cdots \xx_n$ holds for each~$i$. By definition,
  the latter relation always implies $\xx_i \cov(\SS)
  \xx_{i+1}$. By Lemma~\ref{L:cov1}(ii), the converse
  implication is also true: a descending induction on~$j$
  shows that $(\forall i \ge j)(\xx_i \cov(\SS) \xx_{i+1})$
  implies $\xx_i \cov(\SS) \xx_{i+1} \cdots \xx_{p}$.
  Indeed, the conjunction of $\xx_j \cov(\SS) \xx_{j+1}$
  and $\xx_{j+1} \cov(\SS) \xx_{j+2} \cdots \xx_n$
  implies $\xx_j \cov(\SS) \xx_{j+1} \cdots \xx_n$.
\end{proof}

\begin{rema}
  Instead of using $\SS$-simple elements, we could 
  think of simply considering elements of~$\SS$, and 
  constructing a normal form of~$\xx$ starting with
  a maximal divisor of~$\xx$ in~$\SS$. But, then, the
  normal sequences would not necessarily admit the
  local characterization of Prop.~\ref{P:loca}. For
  instance, in the monoid~$M_1$, if we take $\SS =
  \Div(\a^2)$ (a spanning subset that we shall consider
  in Sec.~\ref{S:gars} below), the two sequences
  $(\a, \b)$ and $(\b, \a)$ would be $\SS$-normal, as
  $\a$ is a maximal divisor of~$\a\b$ in~$\SS$,
  and $\b$ is a maximal divisor of~$\a\b$ in~$\SS$,
  but the concatenated sequence $(\a, \b, \a)$ would
  not, as we have $\a^2 \dive \a\b\a$, and therefore
  $\a$ is not a maximal divisor of $\a\b\a$ in~$\SS$.
\end{rema}

We have seen that the normal form of
Prop.~\ref{P:norl} need not be unique in general. We shall
need in Sec.~\ref{S:auto} below the following refinement of
Prop.~\ref{P:norl} that connects the various normal
decompositions of an element:

\begin{prop}\label{P:nopr}
  Assume that $M$ is a \q-atomic cancellative monoid, 
  $\SS$ \q-spans~$M$,
  and $\xx_1$, \pp, $\xx_n$ are $\SS$-simple elements
  of~$M$. Then $\xx_1\cdots \xx_n$ admits a $\SS$-normal
  decomposition $(\xx'_1, \pp, \xx'_{m})$
  such that $m \le n$ holds and, for each~$i$, we
  have
  $\xx_1 \cdots \xx_{f(i)} \dive \xx'_1 \cdots
  \xx'_i$ for some increasing
  mapping $f$ of~$\{1, \pp, m\}$ into~$\{1,$
  \pp, $n\}$ with $f(m) = n$.
\end{prop}

\begin{proof}
  The result is trivial for~$n = 1$. Assume $n = 2$. 
  Applying Prop.~\ref{P:norl} to~$\xx_1\xx_2$, we
  find a $\SS$-normal decomposition of~$\xx_1\xx_2$
  that begins with some multiple~$\xx'_1$ of~$\xx_1$.
  Two cases
  may happen. Either $\xx_1\xx_2$ is $\SS$-simple, and
  $(\xx'_1)$ is the expected decomposition. Or we have
  $\xx'_1 = \xx_1 \yy$ with $\yy \div \xx_2$, and,
  therefore, $\xx_1\xx_2 = \xx'_1
  \xx'_2$ with $\xx_2 = \yy \xx'_2$. By
  Lemma~\ref{L:ridi}, the latter relation forces $\xx'_2$
  to be $\SS$-simple, and, therefore, $(\xx'_1, \xx'_2)$ is
  a $\SS$-normal decomposition of the expected form.
  
  For~$n \ge 3$, we use induction on~$n$. Applying the
  induction hypothesis, we find a $\SS$-normal 
  decomposition
  $(\yy_2, \pp, \yy_p)$ for $\xx_2 \cdots \xx_n$
  and an increasing mapping~$g$ of $\{2, \pp, p\}$
  into $\{2, \pp, n\}$ satisfying
  $\xx_2 \cdots \xx_{g(i)} \dive \yy_2 \cdots \yy_i$
  for $2 \le i \le p$. If $p < n$ holds, we can
  apply the induction hypothesis to~$\xx_1$, $\yy_2$,
  \pp, $\yy_p$, and get the result directly. So, assume
  $p = n$. Then $g$ must be the identity mapping.
  Applying the result with $n = 2$ to $\xx_1 \yy_2$,
  we obtain a $\SS$-normal
  decomposition of length $2$ or $1$. In the latter
  case, we resort to the induction
  hypothesis directly. So, assume that we have
  obtained $(\xx'_1, \yy'_2)$ with $\xx_1 \dive \xx'_1$
  and $\xx'_1 \yy'_2 = \xx_1  \yy_2$. We apply
  the induction hypothesis another time
  to $\yy'_2 \yy_3 \cdots \yy_n$, obtaining a 
  $\SS$-normal decomposition $(\xx'_2, \pp,
  \xx'_{m})$.
  Then $(\xx'_1, \xx'_2, \pp, \xx'_{m})$
  satisfies our requirements. Indeed, by construction, we
  have $\yy_2 \cov(\SS) \yy_3 \cdots \yy_n$, hence, 
  by Lemma~\ref{L:cov1}, $\xx'_1\yy'_2 = 
  \xx_1\yy_2 \cov(\SS) \yy_3
  \cdots \yy_n$, and, as $\xx'_1 \cov(\SS) \yy'_2$ holds
  by construction, $\xx'_1 \cov(\SS) \yy'_2 \yy_3 \cdots
  \yy_n$, hence $\xx'_1 \cov(\SS) \xx'_2$. So the sequence
  $(\xx'_1, \pp, \xx'_{m})$ is $\SS$-normal. The 
  relations $\xx_1 \cdots \xx_{f(i)} \dive \xx'_1 \cdots
  \xx'_i$ follow from the induction hypothesis.
\end{proof}

Although natural, the previous result was not obvious:
putting in normal form a product of two simple
elements might have required say three simple
elements, since the conditions for being normal
discard some decompositions.

We consider now the effect of multiplication on normal
forms, \ie, we try to connect the normal form(s) of an
element~$\xx$ with those of~$\yy\xx$ and~$\xx\yy$,
especially when $\yy$ is $\SS$-simple. As one can expect,
such results will be crucial for constructing an automatic
structure.

\begin{lemm}\label{L:norl}
  Assume that $M$ is a \q-atomic cancellative monoid, 
  and $\SS$ \q-spans~$M$. Let $\xx$,
  $\yy$ be arbitrary elements of~$M$, and $(\xx_1, \pp, 
  \xx_n)$ be a $\SS$-prenormal decomposition of~$\xx$.
  Put $\yy_0 = \yy$, and, inductively, let
  $(\xx'_i, \yy_i)$ be any $\SS$-prenormal decomposition
  of~$\yy_{i-1} \xx_i$. Then $(\xx'_1, \pp, \xx'_n, \yy_n)$
  is a $\SS$-prenormal decomposition of~$\yy\xx$.
  
  If, in addition, $\yy$ is $\SS$-simple and $(\xx_1, \pp, 
  \xx_n)$ is $\SS$-normal, then we may assume that each
  element~$\yy_i$ is $\SS$-simple, and then
  $(\xx'_1, \pp, \xx'_n, \yy_n)$ is a $\SS$-normal
  decomposition of~$\yy\xx$.
\end{lemm}

\begin{figure}
  \includegraphics{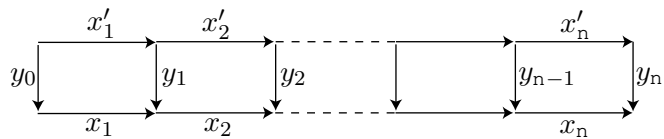}
   \caption{Comparing normal forms
  of~$\xx$ and $\yy\xx$}\label{F:norl}
\end{figure}

\begin{proof}
  (Fig.~\ref{F:norl})
  We have $\xx'_n \cov(\SS) \yy_n$ by construction, so the
  point is to show $\xx'_i \cov(\SS) \xx'_{i+1}$ for each~$i$.
  Assume $\zz \in \SS$ and $\zz \dive \xx'_i \xx'_{i+1}$.
  This implies $\zz \dive \xx'_i \xx'_{i+1}
  \yy_{i+1}$, \ie, $\zz \dive \yy_{i-1} \xx_i \xx_{i+1}$.
  By hypothesis, we have $\xx_i \cov(\SS) \xx_{i+1}$, 
  hence, by Lemma~\ref{L:cov1}(i), $\yy_{i-1}\xx_i
  \cov(\SS) \xx_{i+1}$. So $\zz \dive \yy_{i-1} \xx_i
  \xx_{i+1}$ implies $\zz \dive \yy_{i-1} \xx_i$, \ie,
  $\zz \dive \xx'_i \yy_i$. Now, by hypothesis, we have
  $\xx'_i \cov(\SS) \yy_i$, so we
  deduce $\zz \dive \xx'_i$, hence $\xx'_i
  \cov(\SS) \xx'_{i+1}$, and the sequence $(\xx'_1, 
  \pp, \xx'_n, \yy_n)$ is $\SS$-prenormal.
  
  If $\yy$ and each~$\xx_i$ are $\SS$-simple, we can
  inductively assume that $\yy_{i-1}$ and $\xx'_i$ are
  $\SS$-simple: indeed, in this case,
  Prop.~\ref{P:nopr} guarantees that $\yy_{i-1}
  \xx_i$ admits a $\SS$-normal form of length~2 at most,
  and, if we define $(\xx'_i, \yy_i)$ to be such a $\SS$-normal sequence (with possibly $\yy_i = 1$),
  then induction continues.
\end{proof}

Let us finally consider multiplication on the right. A
similar argument is possible, but it works in the
Gaussian case only.

\begin{lemm}\label{L:norr}
  Assume that $M$ is a Gaussian monoid, 
  $\SS$ \q-spans~$M$. Let $\xx$,
  $\yy$ be arbitrary elements of~$M$, and $(\xx_1, \pp, 
  \xx_n)$ be a $\SS$-prenormal decomposition of~$\xx
  \yy$.
  Put $\yy_n = \yy$, and, inductively, define $\xx'_i $ and
  $\yy_{i-1}$ so that $\yy_{i-1} \xx_i = \xx'_i \yy_i$
  holds and the latter is a left lcm
  of~$\xx_i$ and~$\yy_i$. Then $(\xx'_1, \pp, \xx'_n)$ is
  a $\SS$-prenormal decomposition of~$\xx$.
\end{lemm}

\begin{proof}
  (Fig.~\ref{F:norl} again)
  Let us show that $\xx_i \cov(\SS)
  \xx_{i+1}$ implies $\xx'_i \cov(\SS) \xx'_{i+1}$. Assume
  $\zz \in \SS$ and  $\zz
  \dive \xx'_i \xx'_{i+1}$. As in the proof of
  Lemma~\ref{L:norl}, we deduce $\zz \dive \yy_{i-1}
  \xx_i \xx_{i+1}$, hence,
  as $\xx_i \cov(\SS) \xx_{i+1}$ implies $\yy_{i-1}\xx_i
  \cov(\SS) \xx_{i+1}$, $\zz \dive \yy_{i-1}\xx_i$, \ie,
  $\zz \dive \xx'_i
  \yy_i$. By construction, $\yy_i$ and $\xx'_{i+1}$ have
  no common divisor; {\it in the Gaussian case}, this implies
  that every common divisor of~$\xx'_i \yy_i$ and $\xx'_i
  \xx'_{i+1}$ is a divisor of~$\xx'_i$. So, we have $\zz
  \dive \xx'_i$, and $\xx'_i \cov(\SS) \xx'_{i+1}$.
  
  Finally, we observe that $\yy_0 = 1$ necessarily holds,
  as, by the results of~\cite{Dgk}, $\yy_0 \xx_1 \cdots
  \xx_n$ has to be a left lcm of $\xx_1 \cdots \xx_n$
  and~$\yy$, hence to equal $\xx_1 \cdots \xx_n$.
\end{proof}

\begin{exam}
  When lcm's need not exist, the previous argument fails,
  and so does the result itself. For instance, let us
  consider $M_1$ again. Choose $\xx = \a^3$, $\yy= \b$.
  Then $(\a^2,
  \a)$ is a normal decomposition of~$\xx$. One
  possibility according to 
  Lemma~\ref{L:norr} is to define $\yy_0 = 1$, $\yy_1
  = \b$, $\xx'_1 = \b$, $\xx'_2 = \a$. Indeed, $\a\b=\b\a$
  is a left mcm of~$\a$ and~$\b$, and $1 \a^2 = \b \b$ is
  a left mcm of~$\a^2$ and~$\b$. Now, $(\b, \a)$
  is not a (pre)normal sequence, as we have
  $\a \in \SS$, $\a \dive \b\a$, and $\a \not\dive \b$,
  hence $\b \notcov(\SS) \a$.
\end{exam}

\section{Garside elements}\label{S:gars}

As was recalled above, if $M$ is a thin Gaussian
monoid, \ie, a Garside monoid, then the lcm~$\D$ of
all primitive elements plays an important r\^ole.
Technically, the point is that the left divisors of~$\D$
coincide with its right divisors, which implies in
particular that conjugation by~$\D$ gives an
automorphism of~$M$, and that some power
of~$\D$ belongs to the center of~$M$. Conversely, it
is proved in~\cite{Dgk} that, if $M$ is a Gaussian
monoid and $\D$ is an element of~$M$ such that the
left divisors  of~$\D$ coincide with its right divisors
and they generate~$M$, then these divisors of~$\D$
span~$M$, and, therefore, $M$ is thin, and,
therefore, it is a Garside monoid.

In the general case, there seems to be no reason why the
existence of a finite spanning set should imply
the existence of an element~$\D$ with similar properties.
Even worse, Prop.~\ref{P:node} shows that the existence
of such an element is impossible in the non-Gaussian
case if we require both closure under mcm and left divisors.

However, we shall see now how to define an
appropriate notion of a Garside element which may
exist in the non-Gaussian case, and extends the usual
notion in the Gaussian case. We shall then prove in the
general case a large part of the results established in the
Gaussian case.

\begin{defi}\label{D:gars}
  Assume that $M$ is a cancellative monoid.
  We say that an element~$\D$ of~$M$ is a {\it Garside}
  element if $\Dl(\D)$ is a finite spanning subset
  of~$M$.
\end{defi}

Notice that, if $\D$ is a Garside element in~$M$, then
$M$ must be thin by definition, hence
\q-atomic by Prop.~\ref{P:that}, and every primitive
element of~$M$ must divide~$\D$, since, by
Prop.~\ref{P:micl}, the family~$\PM$ is the least spanning
subset of~$M$, and, therefore, it must be included
in~$\Div(\D)$. Let us mention that most of the
subsequent results could be extended to a {\it
\q-Garside} element, the latter being defined as
an element~$\D$ such that $\Div(\D)$ spans $M$
and is \q-finite.

\begin{lemm}\label{L:auto}
  Assume that $M$ is a thin cancellative monoid,
  and $\D$ is a Garside element in~$M$. Then, for
  every element~$\xx$ in~$\Div(\D)$,
  there exists a unique element~$\xx\dd$ in $\Div(\D)$
  satisfying $\xx \xx\dd = \D$; the mapping $\xx \mapsto
  \xx\dd$ is a permutation of $\Div(\D)$; for $\xx$,
  $\yy \in \Div(\D)$, $\xx$ being a left divisor of~$\yy$
  is equivalent to $\yy\dd$ being a right divisor of~$\xx\dd$.
\end{lemm}

\begin{proof}
  (The argument was already used for
  Prop.~\ref{P:node}.)
  By definition $\xx \dive \D$ means
  that $\xx \xx\dd = \D$ holds for some right
  divisor~$\xx\dd$ of~$\D$, which is unique as $M$ is
  assumed to be (left) cancellative. By hypothesis,
  the family~$\Div(\D)$ spans~$M$
  which contains~$\D$, so, by Prop.~\ref{P:clos},
  it also contains every right divisor of~$\D$,
  so, in particular, $\xx\dd$ belongs
  to~$\Div(\D)$. Then $\xx\dd = \yy\dd$ implies
  $\xx \xx\dd = \D = \yy \yy\dd = \yy \xx\dd$, hence
  $\xx = \yy$, as $M$ is cancellative. This proves
  that $\xx \mapsto \xx\dd$ is an injection
  of~$\Div(\D)$ into itself, hence a bijection as
  $\Div(\D)$ is assumed to be finite. Finally $\yy = 
  \xx \zz$ implies $\xx \xx\dd = \D = \yy \yy\dd = \xx \zz
  \yy\dd$, hence $\xx\dd = \zz \yy\dd$.
\end{proof}

We deduce that, in the Gaussian case, our current
notion of a Garside element coincides with that
considered in~\cite{Dgk}:

\begin{lemm}
  Assume that $M$ is a thin cancellative monoid.

  \noindent (i) If $\D$ is a Garside element
  in~$M$, then the left and the right divisors of~$\D$
  coincide and they generate~$M$. 
  
  \noindent (ii) Conversely, if $M$ is Gaussian
  and $\D$ is an element of~$M$ such that the left
  and the right divisors of~$\D$ coincide and they
  generate~$M$, then $\D$ is a Garside element
  in~$M$.
\end{lemm}

\begin{proof}
  (i) Assume that $\D$ is Garside. By
  Prop.~\ref{P:clos}, every right divisor of~$\D$
  belongs
  to~$\Div(\D)$, hence is a left divisor of~$\D$,
  while, by Lemma~\ref{L:auto}, every element
  of~$\Div(\D)$ belongs to the range of the
  mapping $\xx \mapsto \xx\dd$, hence it is a right
  divisor of~$\D$: so the left and the right
  divisors of~$\D$ coincide.

  (ii) Assume now that $M$ is Gaussian and the left
  and the right divisors of~$\D$
  coincide and they generate~$M$. Assume
  $\xx, \yy \dive \D$ and $\xx \yy'' = \yy \xx''$. Let
  $\xx\yy' = \yy\xx'$ be the lcm of~$\xx$ and~$\yy$. By
  definition, we have $\xx' \dive \xx''$ and $\yy' \dive
  \yy''$. Moreover, $\xx, \yy \dive \D$ implies $\xx\yy'
  \dive \D$. So $\xx\yy'$ and $\yy\xx'$ are left divisors
  of~$\D$, hence they are right divisors of~$\D$ as
  well, and so are $\yy'$ and $\xx'$. Finally, $\xx'$
  and $\yy'$ belong to~$\Div(\D)$, and the latter 
  spans~$M$. So $\D$ is a Garside element.
\end{proof}
  
In the thin Gaussian case, \ie, in a Garside monoid,
there always exists a unique minimal Garside element,
namely the lcm of all primitive elements. In the
general case, we have no such result, but the
following examples show that Garside elements may
still exist.

\begin{exam}
  Consider again the monoid~$M_1$ of
  Example~\ref{X:main}. Let $\D_1 = \a^2$
  and $\D_2 = \a\b$. Then $\D_1$ and $\D_2$ both
  are minimal Garside elements. For instance, we have
  $\Dl(\D_1) = \{1, \a, \b, \a^2\}$, a spanning
  subset of~$M_1$, and the left and the right divisors
  of~$\D_1$ coincide.
  Observe that, in this case, the divisors of~$\D_1$
  properly include the primitive elements.
  
  The reader can check similarly that the monoid~$M_2$
  contains three minimal Garside elements, namely
  $\a^2$, $\a\b$, and $\a\c$, while $M_3$ contains
  one minimal Garside element only, namely~$\b^2$.
\end{exam}

\begin{prop}\label{P:auto}
  Assume that $M$ is a thin cancellative monoid,
  and $\D$ is a Garside element in~$M$. The
  mapping $\xx \mapsto \xx\ddd$
  extends into an automorphism~$\phi_\D$ of~$M$ 
  and we have
  \begin{equation}\label{E:auto}
    \xx \D = \D \phi_\D(\xx)
  \end{equation}
  for every~$\xx$ in~$M$. The
  automorphism~$\phi_\D$
  globally preserves $\Div(\D)$, the units, the atoms, the
  primitive elements, and the simple elements
  of~$M$. The order of~$\phi_\D$ is a finite
  integer~$e$, and the element~$\D^e$ belongs to
  the center of~$M$, which therefore is not trivial.
\end{prop}

\begin{proof}
  By~Lemma~\ref{L:auto}, the mapping $\xx
  \mapsto \xx\ddd$ is a
  permutation of~$\Div(\D)$, and it has a finite order
  say~$e$. By definition, we have
  $\xx \D = \xx (\xx\dd \xx\ddd) = (\xx \xx\dd) \xx\ddd
  = \D \xx\ddd$ for every~$\xx$ in~$\Div(\D)$.
  Assume $\xx_1 \cdots \xx_p = \yy_1 \cdots \yy_q$ with
  $\xx_1$, \pp, $\yy_q \in \Div(\D)$. Using the
previous
  remark, we obtain
  \[
    \D \xx_1\ddd \cdots \xx_p\ddd
    = \xx_1 \cdots \xx_p \D
    = \yy_1 \cdots \yy_q \D
    = \D \yy_1\ddd \cdots \yy_q\ddd,
  \]
  hence $\xx_1\ddd \cdots \xx_p\ddd = \yy_1\ddd
  \cdots \yy_q\ddd$ by cancelling~$\D$.
  Thus putting $\phi_\D(\xx_1 \cdots \xx_p) =
  \xx_1\ddd \cdots \xx_p\ddd$ yields a well defined
  mapping. As $\Div(\D)$ generates~$M$, the
  mapping~$\phi_\D$ is defined everywhere on~$M$, and,
  by construction, it is an endomorphism and
  \eqref{E:auto} is satisfied. Then, $\phi_\D^e$ is also
  an endomorphism, and it is the identity on~$\Div(\D)$,
  so it is the identity everywhere. Hence $\phi_\D$
  must be an automorphism. Moreover, \eqref{E:auto}
  inductively implies $\xx \D^k = \D^k \phi_\D^k(\xx)$
  for every positive~$k$ and every~$\xx$, so, in
  particular, $\xx \D^e = \D^e \xx$ for every~$\xx$, \ie,
  $\D^e$ commutes with every element of~$M$.
  Finally, we apply Prop.~\ref{P:prau}.
\end{proof}

\begin{exam}
  Different Garside elements may give rise to different
  automorphisms. For instance, in $M_2$, the
  automorphism~$\phi_{\a^2}$ is the identity, while
  $\phi_{\a\b}$ and $\phi_{\a\c}$ have
  order~$3$, and they correspond to the cyclic
  permutations $(\a, \c, \b)$ and $(\a, \b, \c
  )$ of the atoms respectively.
\end{exam}

\begin{prop}\label{P:ecmg}
  Assume that $M$ is a thin cancellative monoid,
  and $\D$ is a Garside element in~$M$. Then any two
  elements of~$M$ admit a common multiple; more
  precisely, for $\xx \in \Div(\D)^p$ and
  $\yy \in \Div(\D)^q$, we have $\xx\yy' = \yy\xx'$
  for some $\xx'$ in~$\Div(\D)^p$ and
  $\yy'$ in~$\Div(\D)^q$.
\end{prop}

\begin{proof}
  The proof of Prop.~\ref{P:excm} shows
  that, if $\SS$ spans~$M$ and any two
  elements of~$\SS$ admit a common multiple, then
  two elements~$\xx$ of~$\SS^p$ and $\yy$ of~$\SS^q$
  admit a common multiple $\xx\yy' = \yy\xx'$ with
  $\xx' \in \SS^p$ and $\yy' \in \SS^q$. Here we apply
  the result to the spanning subset~$\Div(\D)$. The only point
  to check is the result in the case $p = q = 1$. Now, for 
  $\xx \dive \D$ and $\yy \dive \D$, we can take
  $\xx' = \yy\dd$ and $\yy' = \xx\dd$.
\end{proof}

For a while let us write $\Dr(\xx)$ for the set of all
right divisors of~$\xx$.

\begin{lemm}\label{L:powe}
  Assume that $M$ is a thin cancellative monoid,
  and $\D$ is a Garside element in~$M$. Then, for every
  positive integer~$k$, we have $\Dl(\D^k) = \Dr(\D^k)
  = \Div(\D)^k$, and, therefore, $\D^k$ is a Garside
  element.
\end{lemm}

\begin{proof}
  We prove the three relations $\Dr(\D^k) \ince
  \Div(\D)^k \ince \Dl(\D^k) \ince \Dr(\D^k)$.
  First, by Prop.~\ref{P:powe}, the set~$\Div(\D)^k$ 
  spans~$M$, and it contains~$\D^k$, so,
  by Prop.~\ref{P:clos}, it also contains every right
  divisor of~$\D^k$.
  
  The second inclusion is proved using induction
  on~$k$. The result is trivial for $k = 1$. Assume $k
  \ge 2$, and let $\xx \in \Div(\D)^k$, say $\xx =
  \xx_1 \xx'$
   with $\xx_1 \dive \D$ and $\xx'
  \in\Div(\D)^{k-1}$ (Fig.	\ref{F:delt}). By
  construction, $\xx_1\dd$ belongs to~$\Div(\D)$, so, by
  Proposition~\ref{P:ecmg}, we have $\xx' \yy =
\xx_1\dd
  \xx''$ for some $\yy \in \Div(\D)$ and
  $\xx'' \in \Div(\D)^{k-1}$. By induction hypothesis,
  we have $\xx'' \dive \D^{k-1}$, and, therefore,
  \[
    \xx \dive \xx \yy = \xx_1 \xx' \yy = \xx_1 \xx_1\dd \xx''
    = \D \xx'' \dive \D \op \D^{k-1} = \D^k.
  \]
  
\begin{figure}
  \includegraphics{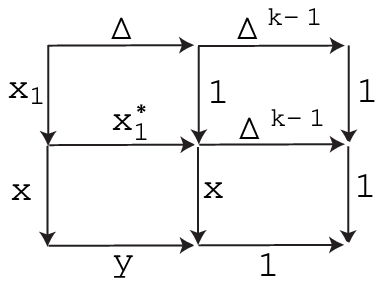}
 \caption{Divisors of $\D^k$}\label{F:delt}
\end{figure}

For the third inclusion, assume $\xx \dive \D^k$, 
  say $\xx \yy = \D^k$. We find
  \[
    \phi_\D^{-k}(\yy) \xx \yy = \phi_\D^{-k}(\yy) \D^k
    = \D^k \phi_\D^k(\phi_\D^{-k}(\yy)) = \D^k \yy,
  \]
  hence $\phi_\D^{-k}(\yy) \xx = \D^k$,  which shows
  that $\xx$ is a right divisor of~$\D^k$.
\end{proof}

\begin{prop}
  Assume that $M$ is a thin cancellative monoid,
  and $\D$ is a Garside element in~$M$. Then any two
  elements of~$M$ admit a common left multiple. 
\end{prop}

\begin{proof}
  Assume $\xx$, $\yy \in M$. Then both $\xx$ and
  $\yy$ belong to~$\Div(\D)^k$ for $k$ large enough.
  By Lemma~\ref{L:powe}, this implies that $\D^k$
  is both a common right multiple, and a common left
  multiple of~$\xx$ and~$\yy$.
\end{proof}

If a thin cancellative monoid~$M$ contains a Garside
element, then, by Proposition~\ref{P:ecmg}, it is a Ore
monoid, and, therefore, it embeds in a thin group of
(right) fractions~$G$. Using the Garside element, we can 
also express every element of~$G$ as a left fraction whose
denominator is a power of~$\D$.

\begin{prop}\label{P:frac}
  Assume that $M$ is a thin cancellative monoid,
  and $\D$ is a Garside element in~$M$. Then $M$
  embeds in a group of fractions~$G$; every element
  of~$G$ admits a unique decomposition $\D^{-k}\xx$
  with $k \in \ZZ$ and $\xx \in M$
  satisfying $\D \not\dive \xx$.
\end{prop}

\begin{proof}
  Let $\zz = \xx \yy\ii$ be
  an element of~$G$. As $\yy$ belongs to~$\Div(\D)^\ell$
  for some positive~$\ell$, we have $\yy \dive \D^\ell$ by
  Lemma~\ref{L:powe}, say $\yy\xx' = \D^\ell$. Then
  we find
  \[
    \zz= \xx \xx' {\xx'}\ii \yy\ii = \xx\xx' \D^{-\ell}
    = \D^{-\ell} \phi_\D^{-\ell}(\xx\xx'),
  \]
  \ie, $\zz = \D^{-\ell} \zz_0$ for some~$\zz_0$
  in~$M$. Assume $p \le \ell$ and $\yy = \D^p\zz
  \in M$. 
  Then, in~$M$, we have $\zz_0 = \D^{\ell - p}
  \yy$, hence $\ell - p \le \n(\zz_0)$. 
  Thus the set $\{p \in \ZZ \, ; \, \D^p \zz  \in M\}$
  must have a least element, say~$k$. Then, by
  construction,
  $\zz$ can be expressed as $\D^{-k}\xx $ for some~$\xx$
  in~$M$. As any relation $\xx = \D\xx'$ in~$M$
  would
  imply $\zz  = \D^{-k+1}\xx'$ and contradict the
  definition of~$k$, we have $\D \not\dive\xx$.
  Finally,
  $\D^{-k}\xx = \D^{-k'}\xx' $ with $k' > k$ implies
  $\D^{-k + k'} \xx = \xx'$ in~$M$, hence $\D 
  \dive \xx'$, showing the uniqueness of the
  decomposition~$\D^{-k} \xx$ when $\D
  \not\dive \xx$ is required.
\end{proof}

If $\D$ is a Garside element in a monoid~$M$, then, by
hypothesis, the set~$\Div(\D)$ spans~$M$, and,
therefore, there exist the associated notions of
a $\Div(\D)$-simple element and a $\Div(\D)$-normal
sequence. For simplicity, we call them $\D$-simple and
$\D$-normal respectively, and we write $\cov(\D)$
for $\cov({\Div(\D)})$. Using
Prop.~\ref{P:frac} and the results of
Sec.~\ref{S:nofo}, we obtain:

\begin{prop}\label{P:nogr}
  Assume that $M$ is a thin cancellative monoid,
  $\D$ is a Garside element in~$M$, and $G$ is the
  group of fractions of~$M$. Then every element of~$G$
  admits a decomposition $\D^{-k} \xx_1 \cdots \xx_p$
  where $k$ is a uniquely determined integer and 
  $(\xx_1, \pp, \xx_p)$ is a $\D$-normal sequence with
  $\xx_1 \not\simeq \D$.
\end{prop}

\begin{proof}
  The only point to establish is that, if $(\xx_1, \pp,
  \xx_p)$ is a $\D$-normal sequence, then $\D \not\dive 
  \xx_1 \cdots \xx_p$ is equivalent to $\xx_1
  \not\simeq \D$.
  The condition is obviously necessary. Conversely, as
  $\D$ belongs to~$\Dl(\D)$, the relation
  $\D \dive \xx_1 \cdots \xx_p$ implies  $\D \dive
  \xx_1$ by definition of a $\D$-normal sequence. Now,
  as $\D$ is divisible by every element of~$\Div(\D)$,
  no proper multiple of~$\D$ may be $\D$-simple,
  and $\D \dive \xx_1$ implies $\D \simeq \xx_1$
  when $\xx_1$ is $\D$-simple.
\end{proof}

\begin{exam}\label{X:desi}
  Even if we use a minimal Garside element, the
  $\D$-normal form need not coincide with the
  ($\PM$)-normal form in general. For instance, consider
  once more the monoid~$M_1$ of Example~\ref{X:main}.
  We have seen that $\D_1 = \a^2$ is a minimal Garside
  element in~$M_1$. Then the $\D_1$-simple elements
  coincide with the simple elements: there are five of them,
  namely $1$, $\a$, $\b$, $\a^2$, and $\a\b$. Now, the
  relations~$\cov()$ and $\cov(\D_1)$ do not
  coincide, because we have $\Div(\D_1) = P_{M_1} \cup
  \{\D_1\}$. It follows that the $\D_1$-simple
elements
  are determined by their divisors in~$\Div(\D_1)$, while
  they are not determined by their primitive divisors: both
  $\a^2$ and $\a\b$ are divisible by~$1$, $\a$, $\b$, but
  only $\a^2$ is divisible by~$\a^2$. As a consequence,
  the $\D_1$-normal form is unique, while we have seen
  the normal form is not.
\end{exam}

\section{Automatic structure}\label{S:auto}

In the Gaussian case, \ie, when lcm exist, thinness
implies the existence of a Garside element, and 
the latter implies the existence of an automatic
structure for the associated group of fractions. We
shall show now that the latter result extends to more
general cases: indeed, we shall prove that, under
suitable hypotheses, the normal form of
Prop.~\ref{P:nogr} is associated with an automatic
structure. 


The first steps, namely proving that the normal
decompositions make a regular language, are easy.

\begin{prop}\label{P:lang}
  Assume that $M$ is a thin cancellative monoid,
  $\D$ is a Garside element in~$M$, and $G$ is the
  group of fractions of~$M$. Let $\SD$ denote the set of
  all $\D$-simple elements in~$M$. Then the language
  consisting of  all
  normal sequences in the sense of
  Prop.~\ref{P:nogr} is regular.
\end{prop}

\begin{proof}
  By Prop.~\ref{P:fisi}, there are finitely many
  $\D$-simple elements, \ie, the set~$\SD$ is finite.
  Put $A = \SD \cup \{\D\ii\}$. 
  A word over~$A$, \ie, a finite sequence
  $(\xx_1, \pp, \xx_n)$ of letters, is a normal form if
  and only if the following requirements are obeyed:
  
  - a letter~$\D\ii$ or~$\D$ cannot follow any other letter;
  
  - a letter~$\xx$ in~$\SD - \{\D\}$ may follow only 
  $\D\ii$ or one of those (finitely many) letters~$\yy$
  in~$\SD$ that satisfy $\yy \cov(\D)
  \xx$.
  
  \noindent Define a state set $Q$ to be $A \cup \{1,
  \bot\}$, where
  $1$ is an initial state and $\bot$ is a failure state, 
  and a transition function $F
  : Q \times A \rightarrow Q$ by
  \begin{center}
    \renewcommand{\arraystretch}{1.10}
     \begin{tabular}{ c || c  | c | c }
       Q $\downarrow \quad A \rightarrow$
       & \quad $\xx \not= \D, \D\ii$ \quad
       & \quad $\D$ \quad 
       & \quad $\D\ii$ \quad \\ \hline\hline
       \quad $\yy \not=  \D, \D\ii$ \quad 
       & $\vcenter{\hsize = 2.2cm
       \strut $\bot$ for $\yy \notcov(\D) \xx$ \par
       \strut $\xx$ for $\yy \cov(\D) \xx$}$
       & $\bot$ & $\bot$ \\ \hline
       $\D$ & $\xx$ & $\D$ & $\bot$ \\ \hline
       $\D\ii$ & $\xx$ & $\bot$ & $\D\ii$ \\ \hline
       $1$ & $\xx$ & $\D$ & $\D\ii$ \\ \hline
       $\bot$ & $\bot$ & $\bot$ & $\bot$
      \end{tabular}
  \end{center}
  Then the finite state automaton $(Q, A, F, 1, Q - \{\bot\})$
  recognizes the language of $\D$-normal forms (see for
  instance \cite{Eps} for definitions).
\end{proof}

Provided all $\D$-normal forms have the same
length, we can readily apply the method of~\cite{Bra}
or~\cite{Xuu}, and deduce:

\begin{coro}
  Assume that $G$ is the group of fractions of a
  cancellative monoid~$M$ that admits a
  Garside element~$\D$ such that all
  $\D$-normal forms of an element have the same length.
  Then $G$ has rational growth, \ie, the number of
  elements of~$G$ with a $\D$-normal form of
  length~$n$ is a rational function of~$n$.
\end{coro}

If $G$ is a group generated by a family~$A$, we denote
by~$\G_A(G)$ the Cayley graph of~$G$ with respect
to~$A$, \ie, the labelled graph whose vertices are the
elements of~$G$ and there exists a $\zz$-labelled edge
from~$\xx$ to~$\yy$ if $\yy = \xx\zz$ holds in~$G$. For
$\xx$, $\yy \in G$, the distance $\dist_{A, G}(\xx,
\yy)$ between~$\xx$ and~$\yy$ in~$\G_A(G)$ is the
minimal length of an unoriented path from~$\xx$
to~$\yy$. 

\begin{defi}
  Assume that $G$ is a group generated by ~$A$. The
  {\it synchronous distance} between two
  words on~$A$, \ie, two sequences 
  of letters in~$A$, say $(\xx_1, \pp, \xx_p)$ and $(\yy_1,
  \pp, \yy_q)$, is defined to be the supremum of
  the numbers
  \[
    \dist_{A, G} (\xx_1 \cdots \xx_{\inf(i, p)}, 
    \yy_1 \cdots \yy_{\inf(i, q)})
  \]
  for $1 \le i \le \sup(p, q)$.
\end{defi}

By the results of~\cite{Eps}, the $\D$-normal form of
Prop.~\ref{P:nogr} is associated with a (left)
automatic structure if and only if the fellow traveller
property (FTP) is satisfied, \ie, for every~$\xx$ in the
group and every~$\yy$ in~$\SD \cup \{\D\ii\}$,

- the synchronous distance between any two
$\D$-normal  decompositions of~$\xx$ is uniformly
bounded, and

- the synchronous distance between a $\D$-normal 
decomposition of~$\xx$ and one of~$\yy\xx$ is
uniformly bounded.

\noindent We shall see that such conditions are satisfied
in good cases. To this end, we shall first establish a
 bound for the distance between the various normal
forms of an element in the monoid. (The notion of the
synchronous distance is extended to the case of the
monoid in the obvious way.)

\begin{lemm}\label{L:anel}
  Assume that $M$ is a \q-atomic cancellative monoid,
  $\SS$ is a \q-spanning subset of~$M$ of
  cardinality~$k$ and, for
  every~$\xx$ in~$M$, the following condition holds:
  \begin{equation}\label{E:anel}
    \text{All $\SS$-normal
     decompositions of~$\xx$ have the same length.}
  \end{equation}
  Then the synchronous distance between any two
  $\SS$-normal  decompositions of an element of~$M$
  is uniformly bounded by~$2(k-1)$.
\end{lemm}

We begin with two auxiliary results.

\begin{lemm}\label{L:copr}
  Assume that $M$ is a (left) cancellative monoid, and $\SS$
  \q-spans~$M$. Then $\xx_1 \cov(\SS) \cdots
  \cov(\SS) \xx_k \cov(\SS) \xx$ implies $\xx_1
  \cdots \xx_k \cov(\SS^k) \xx$.
\end{lemm}

\begin{proof}
  We use induction on~$k \ge 0$. Assume  $\zz \in \SS^k$
  and $\zz \dive \xx_1 \cdots \xx_k \xx$. For $k =
  0$, \ie, for $\zz = 1$, the result is 
  vacuously true. Otherwise, write $\zz  = \zz_1 \zz'$, with 
  $\zz_1 \in \SS$ and $\zz' \in \SS^{k-1}$. By
  Lemma~\ref{L:cov1}(ii), $\xx_1 \cov(\SS) \cdots 
  \cov(\SS) \xx_k \cov(\SS) \xx$ implies $\xx_1
  \cov(\SS) \xx_2 \cdots \xx_k \xx$. By
  hypothesis, we have $\zz_1 \dive \xx_1 \cdots 
  \xx_k \xx$, hence $\zz_1 \dive \xx_1$, say $\xx_1 =
  \zz_1 \xx'_1$. Then, by
  Lemma~\ref{L:cov1}(i), we have $\xx'_1\xx_2 \cov(\SS)
  \xx_3 \cov(\SS) \cdots \cov(\SS) \xx_k \cov(\SS) \xx$,
  and, as $M$ is (left) cancellative, $\zz' \dive (\xx'_1\xx_2)
  \xx_3
  \cdots \xx_k \xx$. By induction hypothesis, this 
  implies $\zz' \dive   (\xx'_1 \xx_2) \xx_3 \cdots \xx_k$,
  hence $\zz = \zz_1 \zz' \dive\zz_1(\xx'_1 \xx_2) \xx_3
  \cdots \xx_k$, \ie, $\zz \dive \xx_1 \cdots \xx_k$.
\end{proof}

\begin{lemm}\label{L:anem}
  Under the hypotheses of Lemma~\ref{L:anel}, if
  $(\xx_1, \pp, \xx_n)$ is a $\SS$-normal decomposition
  for~$\xx$, and $\xx'_1$ is a maximal $\SS$-simple divisor of~$\xx$, then there exist $\xx'_2$, \pp, $\xx'_k$
  such that $(\xx'_1, \pp, \xx'_k, \xx_{k+1}, \pp, \xx_n)$
  is another $\SS$-normal decomposition of~$\xx$.
\end{lemm}

\begin{proof}
  As $\xx'_1$ is $\SS$-simple,  it belongs to~$\SS^k$ by
  Prop.~\ref{P:fisi}. By
  Lemma~\ref{L:copr}, we have $\xx_1 \cdots \xx_k
  \cov(\SS^k) \xx_{k+1} \cdots \xx_n$, so $\xx'_1 \dive
  \xx$ implies $\xx'_1 \dive \xx_1 \cdots \xx_k$, say
  $\xx_1 \cdots \xx_k = \xx'_1 \yy$. Let $(\xx'_2, \pp, 
  \xx'_{k'})$ be a $\SS$-normal decomposition of~$\yy$.
  By hypothesis, we have $\Dl(\xx'_1) \cap \SS =  \Dl(\xx)
  \cap \SS$, hence $\xx'_1 \cov(\SS) \xx'_2
  \cdots \xx'_{k'}$, so $(\xx'_1, \pp, \xx'_{k'})$ is a $\SS$-normal
  sequence, hence another $\SS$-normal decomposition for
  $\xx_1 \cdots \xx_k$. Then Condition~\eqref{E:anel}
  implies $k' = k$. 
  
  Let us now consider the $\SS$-covering relation between
  $\xx'_k$ and~$\xx_{k+1}$. As in the proof of
  Prop.~\ref{P:nopr}, let
  $\xx'$ be a maximal $\SS$-simple divisor of~$\xx'_k
  \xx_{k+1}$ satisfying $\xx'_k \dive \xx'$. Write $\xx'
  = \xx'_k \zz$. Then 
  $\xx'_1 \cdots \xx'_k \zz$ equals
  $\xx'_1 \cdots \xx'_{k-1} \xx'$, so it belongs to~$\SS^k$,
  and, therefore, by Prop.~\ref{P:nopr}ý, it must
  admit at least one normal form of length $k$ at most. On
  the other hand, we have $\xx_k \cov(\SS) \xx_{k+1}$ and
  $\zz \dive \xx_{k+1}$, hence $\xx_k \cov(\SS) \zz$,
  so, if $\zz$ is not invertible, $(\xx_1, \pp, \xx_k, \zz)$ is
  another $\SS$-normal decomposition of $\xx'_1 \cdots
  \xx'_{k-1} \xx'$.
  Condition~\eqref{E:anel} discards this possibility. Hence,
  $\zz$ must be invertible, \ie, we must have $\xx'_k
  \cov(\SS) \xx_{k+1}$. So the sequence $(\xx'_1, \pp,
  \xx'_k, \xx_{k+1})$ is $\SS$-normal, and, trivially, so is
  $(\xx'_1, \pp, \xx'_k, \xx_{k+1}, \pp, \xx_n)$.
\end{proof}

\begin{figure}
  \includegraphics{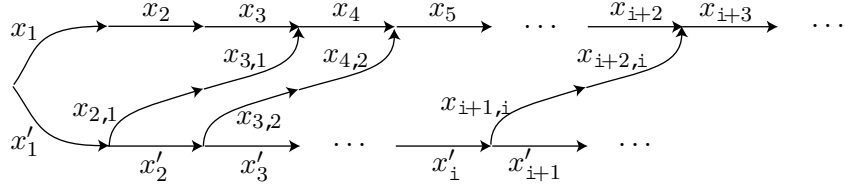}
  \caption{Comparing two normal forms
  of~$\xx$ (here $k = 3$)}\label{F:anel}
\end{figure}

\begin{proof}[Proof of Lemma~\ref{L:anel}]
  Let $(\xx_1, \pp, \xx_n)$ and $(\xx'_1, \pp, \xx'_n)$
  be two $\SS$-normal decomposition of an element~$\xx$
  of~$M$. Applying Lemma~\ref{L:anem} to $(\xx_1,
\pp,
  \xx_n)$ and to~$\xx'_1$, we find $\xx_{2, 1}$,
  \pp, $\xx_{k, 1}$ so that $(\xx'_1, \xx_{2, 1}, \pp,
  \xx_{k, 1}, \xx_{k+1}, \pp, \xx_n)$ is another
  $\SS$-normal decomposition of~$\xx$. Then, applying 
  Lemma~\ref{L:anem} to the latter sequence and
  to~$\xx'_2$, we find $\xx_{3, 2}$, \pp, $\xx_{k+1, 2}$ so
  that $(\xx'_1, \xx'_2, \xx_{3, 2}$,  \pp, $\xx_{k+1, 2},
  \xx_{k+2}, \pp, \xx_n)$ is a $\SS$-normal
  decomposition
  of~$\xx$. Similarly, having found a $\SS$-normal form
  $(\xx'_1, \pp, \xx'_i, \xx_{i+1, i}$,
  \pp, $\xx_{i+k-1, i}$, 
  $\xx_{i+k}, \pp, \xx_n)$ for~$\xx$, applying
  Lemma~\ref{L:anem} to this sequence and
  to~$\xx'_{i+1}$
  yields a new $\SS$-normal decomposition $(\xx'_1, \pp,
  \xx'_{i+1}, \xx_{i+2, i+1}$, \pp, $\xx_{i+k, i+1}, 
  \xx_{i+k+1}, \pp, \xx_n)$. Now, we read on
  Fig.~\ref{F:anel} that, for each~$i$, the distance
  between $\xx_1 \cdots \xx_i$ and $\xx'_1 \cdots \xx'_i$
  is bounded above by~$2(k-1)$, as $\xx_1 \cdots 
  \xx_{i + k-1}$ is a common multiple of these
  elements.
\end{proof}

Applying the previous result to the case of monoids
with a Garside element, we deduce:

\begin{prop}\label{P:vanf}
  Assume that $M$ is a thin cancellative monoid,
  $\D$ is a Garside element in~$M$ with $k$~divisors,
  and $G$ is the group of fractions of~$M$. Assume 
  moreover that, for
  every~$\xx$ in~$M$, the following condition holds:
  \begin{equation}\label{E:vanf}
    \text{All $\D$-normal
     decompositions of~$\xx$ have the same length.}
  \end{equation}
  Then the synchronous distance between any two
  $\D$-normal decompositions of an element of~$G$
  is uniformly bounded by~$2(k-1)$.
\end{prop}

\begin{proof}
  We have seen that, if $(\D\ii, \pp, \D\ii, \xx_1, \pp,
  \xx_n)$, $k$~times~$\D\ii$, and $(\D\ii, \pp, \D\ii,
  \xx'_1, \pp, \xx'_{n'})$, $k'$~times~$\D\ii$, are two
  $\D$-normal decompositions for some element~$\zz$
  of~$G$, then, necessarily, $k = k'$ holds, and, therefore,
  $(\xx_1, \pp, \xx_n)$ and $(\xx'_1, \pp, \xx'_{n'})$
  are two $\D$-normal decompositions for some
  element
  of~$M$. Then we apply Lemma~\ref{L:anel} with
  $\SS = \Div(\D)$.
\end{proof}

The case of left multiplication in the monoid has
already been treated in Lem\-ma~\ref{L:norl}, which gives:

\begin{lemm}\label{L:munf}
  Assume that $M$ is a \q-atomic cancellative monoid,
  $\SS$ \q-spans~$M$, and $\yy$
  is a $\SS$-simple element of~$M$. Then, for every
  element~$\xx$ of~$M$, and every $\SS$-normal
  decomposition of~$\xx$, there exists a $\SS$-normal
  decomposition of~$\yy\xx$  at synchronous distance at
  most~$1$.
\end{lemm}

It remains to extend the result to the group of fractions.

\begin{prop}\label{P:munf}
  Assume that $M$ is a thin cancellative monoid,
  $\D$ is a Garside element in~$M$ with $k$~divisors,
  $G$ is the group of fractions of~$M$, and $\yy$ is a
  $\D$-simple element of~$M$. Then, for every
  element~$\zz$ of~$G$, and every $\D$-normal
  decomposition of~$\zz$, there exists a $\D$-normal
  decomposition of~$\yy\zz$  at synchronous distance at
  most~$3k$.
\end{prop}

\begin{proof}
  Assume first $\yy \in \Div(\D)$.
  Assume $\zz = \D^{-k} \xx$, with $\xx \in M$ and $\D
  \not\dive \xx$. Then we have $\yy\zz= \D^{-k} \yy'
   \xx$ with $\yy' = \phi_\D^{-k}(\yy)$. By
  Prop.~\ref{P:auto}, we have $\yy' \dive \D$, so, in
  particular, $\yy'$ is
  $\D$-simple, and we can apply Lemma~\ref{L:munf} to
  $\yy'$ and any $\D$-normal decomposition $(\xx_1,
  \pp, \xx_n)$ of~$\xx$ to obtain a $\D$-normal
  decomposition $(\xx'_1, \pp, \xx'_n, \yy'_n)$
  of~$\yy'\xx$. There remains one point to check: if
  it contains at least one~$\D\ii$, 
  the sequence $(\D\ii, \pp, \D\ii, \xx'_1$, \pp, $\xx'_n,
  \yy'_n)$ is $\D$-normal only if $\xx'_1$ is not~$\D$: if
  $\xx'_1 = \D$ holds, we must cancel $\xx'_1$ with the
  last~$\D\ii$, and repeat the reduction until we possibly
  find $\xx'_i \not= \D$. As each such reduction increases
  the synchronous distance by~$2$, there could be a
  problem here. Actually, we shall prove that $\xx'_1
  \simeq \xx'_2 \simeq \D$ implies $\D \dive
  \xx_1$, hence $\xx_1 \simeq \D$. Here we use the
  hypothesis that $\yy'$ is not only $\D$-simple, but
  also it is a divisor of~$\D$. First, $\xx'_1 \simeq
  \xx'_2 \simeq \D$ implies $\xx'_1 \xx'_2
  \simeq \D^2$. Indeed, for $\uu \in \UM$,
  we have $\uu \D = \xx \vv$ for some $\xx$
  and~$\vv$ satisfying $\xx \in \Div(\D)$ and $\vv
  \in \UM$, and $\n(\xx) = \n(\D)$ implies $\xx
  \simeq \D$. So we deduce $\D^2 \dive  \yy' \xx_1
  \xx_2$, \ie,  $\yy' {\yy'}\dd \D \dive \yy'
  \xx_1 \xx_2$, hence ${\yy'}\dd \D \dive \xx_1 \xx_2$,
  \ie, $\D \phi_\D({\yy'}\dd) \dive \xx_1 \xx_2$
  which implies $\D \dive \xx_1 \xx_2$, and, finally,
  $\D \dive \xx_1$ as $\xx_1 \cov(\D) \xx_2$ holds
  by hypothesis. So, at most one reduction $\D\ii \op \D$
  may occur, and the synchronous distance between
  the $\D$-normal form of~$\xx$ and that of~$\yy \xx$
  is at most~$3$.
  
  The result for an arbitrary $\D$-simple element~$\yy$
  follows, as, by Prop.~\ref{P:fisi}, any such element
  is the product of at most $k$~elements
  of~$\Div(\D)$.
\end{proof}

Putting Propositions~\ref{P:lang}, \ref{P:vanf},
and~\ref{P:munf} together, we deduce

\begin{prop}\label{P:aust}
  Assume that $G$ is the group of fractions of a
  cancellative monoid~$M$ that admits a
  Garside element~$\D$ such that all
  $\D$-normal forms of an element have the same length.
  Then $G$ is an automatic group.
\end{prop}

The previous result applies in particular to every
thin Gaussian group, \ie, to every Garside
group---hence in particular to every spherical Artin
group. But non-Gaussian groups are also eligible:

\begin{exam}
  Consider once more the groups~$G_1$ and $G_3$
  of Example~\ref{X:maio}. We have seen in
  Example~\ref{X:desi} that the monoid~$M_1$
  contains a Garside element~$\D$ such that the
  $\D$-simple elements are determined by their
  divisors in~$\Div(\D)$. So the associated
  $\D$-normal form is unique, and, therefore, the
  length requirement is satisfied. The argument is
  similar for~$M_3$. So the groups~$G_1$
  and~$G_3$ are automatic.
\end{exam}

The case of~$G_2$ is slightly different. Indeed, in the
monoid~$M_2$, $\a^2$ is a Garside element, but
$\a\b$ and $\a\c$ are
$\a^2$-simple elements with the same divisors
in~$\Div(\a^2)$, namely $1$, $\a$, $\b$, $\c$.
Now, we have the following sufficient condition:

\begin{prop}
  Assume that $M$ is a thin cancellative monoid
  with no nontrivial unit,
  $\D$ is a Garside element in~$M$, and the
  following condition holds in~$M$: If $\xx$ and
  $\xx'$ are distinct $\D$-simple elements
  with the same divisors in~$\Div(\D)$, then
  every common multiple of~$\xx$ and~$\xx'$
  is a multiple of some $\D$-simple common
  multiple of~$\xx$ and~$\xx'$. Then the
  $\D$-normal form is unique,
  and, therefore, the group of fractions of~$M$
  is automatic.
\end{prop}

\begin{proof}
  It suffices to show that, for every~$\xx$ in~$M$,
  there exists a unique $\D$-simple
  element~$\xx_1$ with the same divisors as~$\xx$
  in~$\Div(\D)$. Now, assume that $\xx_1$ and
  $\xx'_1$ satisfy these conditions and are distinct.
  Then, by
  hypothesis, there exists a $\D$-simple
  element~$\xx''_1$ satisfying $\xx_1 \dive \xx''_1
  \dive \xx$ and $\xx'_1 \dive \xx''_1 \dive \xx$.
  Then $\xx'_1 = \xx''_1$ would imply $\xx_1 \div
  \xx'_1$, contradicting the $\D$-simplicity
  of~$\xx'_1$. So we must have $\xx_1 \div
  \xx''_1$, and, therefore, $\Div(\xx_1) \cap
  \Div(\D) \not= \Div(\xx'_1) \cap \Div(\D)
  \ince \Div(\xx) \cap \Div(\D)$, which contradicts
  $\Div(\xx_1) \cap \Div(\D) = \Div(\xx) \cap
  \Div(\D)$. 
\end{proof}

\begin{exam}
  The previous criterion applies to the
  monoid~$M_2$: indeed, for $\D_1 = \a^2$, the
  only problem with~$\D_1$ occurs with the
  $\D_1$-simple
  elements $\a\b$ and $\a\c$. Now, every common
  multiple of~$\a\b$ and~$\a\c$ is a multiple
  of~$\a^2$, \ie, of~$\D_1$. We deduce that 
  $G_2$ is automatic.
\end{exam}

Let us conclude with some open questions.

\begin{ques}
  If $\D$ is a Garside element in a thin cancellative
  mon\-oid~$M$, do all $\D$-normal decompositions of a
  given element of~$M$ necessarily have the same length,
  \ie, is the additional assumption of
  Prop.~\ref{P:aust} superfluous?
\end{ques}

In the Gaussian case, Lemma~\ref{L:norr} gives a
uniform bound for the synchronous distance between
the normal form of~$\xx$ and that of~$\xx\yy$
when $\yy$ is simple. It is then easy to deduce that
the $\D$-normal form of Prop.~\ref{P:aust} gives rise
to a bi-automatic structure---alternatively, we can also
replace in this case the dissymmetric form
$\D^{-k}\xx_1 \cdots \yy_n$ with a symmetric one
$\yy_p\ii \cdots \yy_1\ii \xx_1 \cdots \xx_n$
\cite{Thu, Dgk}. In the general case, the argument
fails, the behaviour of $\D$-normal form with
respect to right multiplication remains unknown, and
so does the existence of an automatic structure
involving a symmetric fractionary decomposition
(defining the latter in the non-Gaussian case seems to
require a uniform bound for the distance between
the possible various mcm's of two elements in
the monoid).

\begin{ques}
  Under the hypotheses of Prop.~\ref{P:aust}, is the
  group~$G$ bi-auto\-matic?
\end{ques}

(In the case of the groups~$G_1$, $G_2$, $G_3$
of Example~\ref{X:maio}, a simple specific argument
gives a positive answer.)

By Prop.~\ref{P:ecmg}, common multiples must exist
in every thin cancellative monoid admitting a Garside
element. In the Gaussian case, \ie, when we assume
not only that common multiples exist, but even that
least common multiples exist, then the lcm of all
primitive elements is a Garside element.

\begin{ques}
  Does every thin cancellative monoid admitting
  common multiples contain a Garside element? More
  precisely, need every mcm of the primitive elements
  be a Garside element?
\end{ques}

Finally, let us mention an open problem dealing with the
Gaussian case:

\begin{ques}
  Is every finitely generated Gaussian group thin, \ie,
  is every finitely generated Gaussian group
  necessarily a Garside group?
\end{ques}

\end{document}